%% file: main.tex
\documentclass[a4paper,10pt]{amsart}

\usepackage{amsrefs}
\usepackage{amsfonts,mathrsfs}
\usepackage{amsthm}
\usepackage{amssymb}
\usepackage{enumerate}
\usepackage{tikz-cd}
\usepackage{cleveref}
\usepackage{enumitem}
\usepackage[obeyFinal,textsize=footnotesize]{todonotes} 
\usepackage{soul}

\theoremstyle{definition}
\newtheorem{Def}{Definition}[section]
\newtheorem{ex}[Def]{Example}
\newtheorem{rem}[Def]{Remark}
\theoremstyle{plain}
\newtheorem{prop}[Def]{Proposition}
\newtheorem{thm}[Def]{Theorem}
\newtheorem*{thm*}{Theorem}
\newtheorem{lem}[Def]{Lemma}

\newtheorem{cor}[Def]{Corollary}
\newtheorem*{cor*}{Corollary}

\newtheorem*{con*}{Conjecture}

\newtheorem*{verm*}{Vermutung}

\theoremstyle{definition}
\newtheorem{assumption}{Assumption}

\newcommand{\supp}{\mathrm{supp}}

\newcommand{\varietyname}{essentially strict Lee--Yang}

\newcommand{\quash}[1]{}

\input{operators}

\input{fonts}

\providecommand{\spann}{\operatorname{span}}
\providecommand{\newt}{\operatorname{Newt}}
\providecommand{\ini}{\operatorname{in}}
\title[Uniformly discrete crystalline counting measures]{A characterization of uniformly discrete crystalline measures with unit weights}
\author{Lior Alon}
\address{Ben-Gurion University of the Negev, Israel}
\email{lioral@bgu.ac.il}
\author{Mario Kummer}
\address{Technische Universit\"at Dresden, Germany} 
\email{mario.kummer@tu-dresden.de}

\begin{document}
	\begin{abstract}

We prove that every uniformly discrete set $\emptyset\ne\Lambda\subset\R^d$ admitting a Fourier summation formula with a discrete set of frequencies $S(\Lambda)\subset\R^d$ arises from the Alon--Kummer--Kurasov--Vinzant construction using a Lee--Yang variety. In particular, every crystalline measure with uniformly discrete support and unit weights is a Fourier quasicrystal coming from the above construction.
It also implies that every set $\Lambda\subset\R^d\subset\C^d$ as above is the common zero set in $\C^d$ of finitely many trigonometric polynomials, extending the one-dimensional theorem of Olevskii and Ulanovskii. As part of the proof, we generalize Meyer's notion of lighthouse measures to higher codimension and introduce lighthouse manifolds, which we characterize, up to a monomial change of coordinates, as the torus parts of Lee--Yang varieties.
\end{abstract}

	\maketitle
	
\section{Introduction}
\subsection*{Notation and Fourier summation formulas}
For a function $f:\R^d\to\C$, we use the convention that its Fourier transform is  
	\begin{equation*}
		\hat{f}(y) = \int_{\mathbb{R}^d} f(x)e^{- 2\pi i \langle x,y\rangle}\,dx.
	\end{equation*}
Poisson's summation formula expresses a fundamental relation between a lattice and its dual: if $\mathcal{L}\subset\R^d$ is a full-rank lattice of covolume one and $\widehat{\mathcal{L}}$ is its dual lattice, then
\[
    \sum_{x\in\mathcal{L}}\hat{f}(x)
    =
    \sum_{k\in\widehat{\mathcal{L}}}f(k)
    \qquad \forall f\in\mathcal{S}(\R^d),
\]
or, equivalently,
\[
    \widehat{\sum_{x\in\mathcal{L}}\delta_x}
    =
    \sum_{k\in\widehat{\mathcal{L}}}\delta_k.
\]
Here $\mathcal{S}(\R^d)$ is the space of Schwartz functions, and the Fourier transform of a tempered measure $\mu$ is defined via $\hat{\mu
}(f):=\mu(\hat{f})$ for all $f\in \mathcal{S}(\R^d)$.

From the Poisson summation formula it is easy to derive a summation formula for any discrete periodic subset of $\R^d$.
A natural question is whether analogous summation formulas exist beyond the periodic setting.

If a given set $\Lambda\subset\R^d$ is uniformly discrete, then its counting measure $$\mu_{\Lambda}:=\sum_{x\in\Lambda}\delta_x$$ is tempered, and as such the Fourier transform $\widehat{\mu_{\Lambda}}$ is a well-defined tempered distribution. In general, it is not a measure. Even when $\widehat{\mu_{\Lambda}}$ is a measure, it is, in general, not a pure-point measure. Recall that, when it is a pure-point measure, we can write  $\widehat{\mu_{\Lambda}}=\sum_{\xi\in S(\Lambda)}c_{\xi}\delta_{\xi}$ with $c_{\xi}\in\C\setminus\{0\}$. In this case we call $S(\Lambda)$ the \emph{spectrum}, and ${c_{\xi}}$ the \emph{Fourier coefficients}. Notice that the \emph{Fourier support} of $\mu_\Lambda$ is the closure $\overline{S(\Lambda)}$. The spectrum and Fourier support agree in the rare case when $S(\Lambda)$ is discrete (i.e., a locally finite set), which is exactly the case considered in this paper. 
\subsection*{Quasicrystals, diffraction, and Lagarias' conjecture}
The question of which sets $\Lambda$ have pure-point $\widehat{\mu_{\Lambda}}$ is intimately related to the question of pure-point diffraction, which gained interest after the discovery of the first quasicrystal.

A Meyer model set $\Lambda\subset\R^d$, also called
a cut-and-project set, is obtained from a suitable
lattice $\mathcal L\subset\R^d\times\R^n$ by selecting
the lattice points whose $\R^n$ coordinates lie in a
prescribed bounded window and projecting them onto
$\R^d$. Let $S$ denote the projection of the dual
lattice $\widehat{\mathcal L}$ onto $\R^d$.
This construction can produce nonperiodic sets
with formal Fourier expansions
\[
    \widehat{\mu_\Lambda}
    \sim
    \sum_{\xi\in S}c_\xi\delta_\xi.
\]
The frequency set $S$ is countable and typically dense
in $\R^d$. The coefficients depend on the window and
may vanish. They need not be locally absolutely
summable, so this formal expansion need not define
a pure-point measure. However, under mild assumptions on the window, the $|c_{\xi}|^2$ are locally summable, in which case $\Lambda$ has a pure-point diffraction measure\footnote{The diffraction measure is defined as the Fourier transform of the autocorrelation measure of $\Lambda$, which exists in this case.} \[
    \gamma_\Lambda
    =
    \sum_{\xi\in S}|c_\xi|^2\delta_\xi.
\]
Physically, when a plane wave $e^{2\pi i\langle x,k\rangle}$ scatters from atoms at $\Lambda$, the diffraction measure describes the intensities of the scattered waves recorded on a distant screen, expressed in the frequency coordinate $\xi=k_{\mathrm{out}}-k$.

These connections between cut-and-project sets, Poisson summation, and diffraction became central to the mathematical theory of quasicrystals \cite{Meyer1995quasicrystals}.

Shechtman's discovery of physical  quasicrystals \cite{shechtman1984metallic}, nonperiodic materials exhibiting sharp diffraction peaks, gave these ideas a physical significance and motivated Bombieri and Taylor's question ``Which distributions of matter diffract?'' \cite{bombieritaylor1986}. Subsequent work of Lagarias, Moody, and others developed the connections between Meyer's theory, aperiodic order, and diffraction \cites{lagarias1996meyer,lagarias2000mathematical,moody2000model}. Lagarias raised a rigidity conjecture, which we state in a simplified version (without weights) as follows: if $\Lambda$ is uniformly discrete, and  $\widehat{\mu_{\Lambda}}$ is a pure-point measure with $S(\Lambda)$ also uniformly discrete, then $\Lambda$ must be periodic \cite{lagarias2000mathematical}. Lev and Olevskii proved Lagarias' conjecture\footnote{They proved the general weighted version in dimension one and, assuming general positive weights, in all dimensions.}\cites{lev2015quasicrystals,lev2017fourier}.
\subsection*{Fourier quasicrystlas and Lee--Yang constructions}

A big open question that was left was whether, in order for the above rigidity result to hold, it is necessary to ask that $S(\Lambda)$ is uniformly discrete, rather than just discrete (i.e., locally finite). Surprisingly, the answer is yes. Kurasov and Sarnak \cite{kurasovsarnak} constructed the first non-periodic uniformly discrete sets $\Lambda\subset\R$ for which $\widehat{\mu_{\Lambda}}$ is a pure point measure with discrete Fourier support $S(\Lambda)$. {The inputs to their construction are so-called \emph{Lee--Yang polynomials}. These are multivariate complex polynomials that satisfy a certain stability condition, i.e., they do not vanish on certain regions of $\C^n$.} Furthermore, the works of Olevskii and Ulanovskii \cite{olevskii2020fourier} and the first author with Cohen and Vinzant \cite{AlonCohenVinzant} proved that the Kurasov--Sarnak construction is in fact the only way to construct discrete sets with discrete $S(\Lambda)$ in dimension one, assuming an additional temperateness condition. 

The authors, together with Kurasov and Vinzant, introduced a general construction of uniformly discrete sets $\Lambda\subset\R^d$, in all dimensions $d$, for which the Fourier support $S(\Lambda)$ is discrete \cite{alon2025higher}. {In this construction, Lee--Yang polynomials are replaced by \emph{strict Lee--Yang varieties}. These are algebraic varieties that avoid certain regions in $\C^n$ and which serve as a higher-codimension analogue of Lee--Yang polynomials.} 

\subsection*{Lighthouse measures}
The Kurasov--Sarnak construction led Meyer to identify another form of simultaneous localization, now on the Fourier-dual pair $(\T^n,\mathbb Z^n)$: a positive singular measure on the torus whose Fourier support lies in a proper double cone, which he called a \emph{lighthouse measure} \cite{Meyer2023}. The elementary example is when $\nu$ is the Haar measure of a codimension-one subtorus $H\subset\T^n$, in which case $\supp(\widehat\nu)$ is the rank-one subgroup of $\Z^n$, defined by $k\in\Z^n$ such that $ z^k=1$ for all $ z\in H$. The notation $z^k$ stands for $\prod_{j=1}^n z_j^{k_j}$. 

In general, let $\nu$ be a lighthouse measure supported on
a hypersurface $M\subset\T^n$. A direction $\ell\in\R^n$
that is transverse to $M$ everywhere, with $\Q$-linearly independent entries,
leads, through the intersection-projection mechanism, to
a formal weighted summation formula
\[
    \widehat{\sum_{t\in\Lambda}a_t\delta_t}
    \sim
    \sum_{k\in\Z^n}\widehat\nu(k)
    \delta_{\langle k,\ell\rangle},
    \qquad
    \Lambda=\{t\in\R\mid\exp(2\pi it\ell)\in M\}.
\]
Here $\exp(z):=(e^{z_1},\ldots,e^{z_n})$, and the weights
$a_t$ are determined by $\nu$. For a certain choice of $\nu$, the
$\ell$-directed measure on $M$, all weights equal one,
giving the counting measure $\mu_\Lambda$.
Our earlier work \cite{alonkummer} shows {that, under mild assumptions, such hypersurfaces are necessarily algebraic, and, up to a choice of coordinates, the zero set of a Lee--Yang polynomial.}

\subsection*{The present work}
We prove that every uniformly discrete set
$\Lambda\subset\R^d$ whose counting measure is crystalline (i.e., $\widehat{\mu_{\Lambda}}$ is a pure-point measure with discrete support)
arises from the Alon--Kummer--Kurasov--Vinzant construction
using a strict Lee--Yang variety. In particular, every
such counting measure is a Fourier quasicrystal:
the additional temperateness condition is automatic.
In the course of the proof, we introduce \emph{lighthouse manifolds}, generalizing
lighthouse hypersurfaces to submanifolds of higher
codimension. We show that every such $\Lambda$ arises
from a lighthouse manifold and that every lighthouse
manifold is the torus part of an algebraic variety
of Lee--Yang type.

\subsection{Preliminaries}
Our hypothesis throughout the paper is the following.
\begin{assumption}[Crystalline counting measure]\label{ass:crystalline}
The set $\Lambda\subset\R^d$ is non-empty, and $\delta$-uniformly discrete:  $\|x-x'\|\ge\delta>0$ for all distinct $x,x'\in\Lambda$. We also assume that $\widehat{\mu_{\Lambda}}$ is a pure-point measure
\begin{equation}\label{eq: Fourier measure}
    \widehat{\mu_\Lambda}=\sum_{\xi\in S(\Lambda)}c_\xi\delta_\xi,
    \qquad c_\xi\in\C\setminus\{0\},
\end{equation}
with discrete (i.e., locally finite) support $S(\Lambda)\subset\R^d$.
\end{assumption}
Uniform discreteness makes $\mu_\Lambda$ tempered, so its distributional Fourier transform is well-defined. A tempered measure $\mu$ with locally finite support whose Fourier transform $\widehat{\mu}$ is a measure with locally finite support is called \emph{crystalline} \cite{meyerquestion}. It is a \emph{Fourier quasicrystal} (FQ) if, in addition, also $|\mu|$ and $|\widehat{\mu}|$ are tempered \cites{lev2015quasicrystals,lev2017fourier}. For a crystalline counting measure $\mu_{\Lambda}$, the additional requirement for being a Fourier quasicrystal is 
\begin{equation}\label{ass:growth}
    \sum_{\xi\in S(\Lambda)\cap B_R}|c_\xi|
    \le CR^m \qquad  \forall R>1,
\end{equation}
for some positive constants $C,m>0$. We also call a set $\Lambda$ an FQ when $\mu_{\Lambda}$ is a Fourier quasicrystal.
Gon\c{c}alves \cite{goncalves2023classification} showed that in dimension one, the polynomial bound in  \eqref{ass:growth} can be replaced with a weaker exponential growth bound, in a more general setting allowing positive weights that have a constant positive lower bound.

\begin{Def}
	A set $\Lambda \subset \R^d$ is \emph{Delone} if it is uniformly discrete and relatively dense. Recall that $\Lambda$ is \emph{relatively dense} if there exists $R > 0$ such that $
		\bigcup_{x \in \Lambda} B_R(x) = \R^d$, namely the balls of radius $R$ around points of $\Lambda$ cover $\R^d$. 
\end{Def}
 Delone Fourier quasicrystals $\Lambda\subset\R^d$ have been constructed in \cite{kurasovsarnak} for $d=1$ and in \cites{alon2025higher, lawtontsikh} for arbitrary $d$. 
 The input to the construction in \cite{alon2025higher} is:{
\begin{itemize}
	\item a subset $X\subset(\C^*)^n$, and 
	\item a real matrix $L\in\R^{n\times d} $.
\end{itemize}
Here $\C^*=\C\smallsetminus\{0\}$, and $\T\subset\C$ is the unit circle.} { The map $z\mapsto \exp(z)$ is the coordinate-wise exponentiation sending $\C^n$ to $(\C^*)^n$, and in particular, the map
$x\mapsto \exp(2\pi iLx)$ sends $\R^d$ into $\T^n$.}
We consider the preimage of $X$ under the latter map:
\begin{equation} 
	\Lambda(X,L): =\{x\in\R^d\mid \exp(2\pi iLx)\in X\}.
\end{equation}
In order for $\Lambda(X,L)$ to be a Delone Fourier quasicrystal it suffices that all $d\times d$ minors of $L$ are positive and $X$ is a strict Lee--Yang variety as defined below.

\begin{Def}[Strict Lee--Yang variety]\label{def:strictLY}
	An equidimensional closed algebraic subvariety $X$ in $ (\C^*)^n $, of dimension $ c=n-d\ge 0 $, is called a \emph{strict Lee--Yang variety} if it satisfies the following three conditions:
	\begin{enumerate}
		\item $X$ is invariant under the coordinate-wise involution $z\mapsto 1/\overline{z}$;
        \item $X(\T):=X\cap\T^n$ is contained in the smooth part of $X$;
        \item { the \emph{amoeba of $X$},
        \begin{equation*}
            \cA(X):=\{\log|z|\in\R^n\mid z\in X\},\quad \text{with}\  \log|z|:=(\log|z_1|,\ldots,\log|z_n|),
        \end{equation*}
        intersects the range $L\R^d$ of $L$ only at the origin, $$\cA(X)\cap L\R^d=\{0\},$$
        for every full-rank matrix $L\in\R^{n\times d}$ with only non-negative $d\times d$ minors.
        }
	\end{enumerate}
\end{Def}
 The Alon--Kummer--Kurasov--Vinzant construction states:
\begin{thm}[\cite{alon2025higher}*{Theorem 2.3}]\label{thm:akkvmain}
	Let $X\subset(\C^*)^n$ be a strict Lee--Yang variety of dimension $c=n-d$ and let $L\in\R^{n\times d}$ be a real matrix all of whose $d\times d$ minors are positive. Then the set $\Lambda=\Lambda(X,L)$ is  Delone and it is a Fourier quasicrystal. 
\end{thm}
 \section{Main results}

We prove the converse of Theorem \ref{thm:akkvmain}, under the weaker assumption that the counting measure of $\Lambda$ is only crystalline.
\begin{thm}\label{thm:main}
If $\Lambda\subset\R^{d}$ satisfies Assumption~\ref{ass:crystalline}, then $$\rank(\spann_{\Z}(S(\Lambda)))=n<\infty,$$ and there exists a strict Lee--Yang variety $X\subset(\C^*)^n$ of dimension $c=n-d$ and a real matrix $L\in\R^{n\times d}$ all of whose $d\times d$ minors are positive such that
    \begin{equation} 
	\Lambda=\Lambda(X,L): =\{x\in\R^d\mid \exp(2\pi iLx)\in X\}.
\end{equation} In particular, $\Lambda$ is Delone and its counting measure is a Fourier quasicrystal.
\end{thm}
For Delone Fourier quasicrystals, the following consequence was conjectured in \cite{lawtontsikh}*{Page 2}.
\begin{cor}
    Every set satisfying Assumption~\ref{ass:crystalline} is the common zero set of finitely many trigonometric polynomials:  if $p_1,\ldots,p_m$, $m\ge d$,  are polynomials whose common zero set is $X$, then $f_{j}(x)=p_{j}(\exp(2\pi i Lx))$ are  trigonometric polynomials whose common zeros in $\C^d$ are exactly $\Lambda$, and in particular, real.
\end{cor}
\subsection{Outline of the proof}
Let $\Lambda$ satisfy Assumption~\ref{ass:crystalline}, and put $\mu=\mu_\Lambda$.
\subsubsection{The spectrum lies in a finitely generated group. }
We first prove in Section \ref{sec:toraltype} that the spectrum of $\mu$ is contained in a finitely generated subgroup of $\R^d$. The proof relies on an interpolation result from \cite{lev2017fourier} that essentially goes back to Ingham and Kahane \cites{ingham,kahane}. 
\subsubsection{Lawton's torus construction}
{Let $\Gamma=\spann_{\Z}(\supp\widehat\mu)$, let $n=\rank(\Gamma)$, and choose $L\in\R^{n\times d}$ whose rows form a $\Z$-basis of $\Gamma$. Then
\[
    L^t\colon\Z^n\longrightarrow\Gamma
\]
is a group isomorphism, and the homomorphism $\varphi\colon\R^d\to\T^n$,
$\varphi(x)=\exp(2\pi iLx)$, defines an ergodic $\R^d$-action on $\T^n$.
Using Lawton's construction \cite{lawton-inverse}, we associate to $\Lambda$ a submanifold $M\subset\T^n$ of complementary dimension $c=n-d$, together with a measure $\nu$ on $M$, such that
\[
    M=\overline{\varphi(\Lambda)},
    \qquad
    \Lambda=\varphi^{-1}(M),
    \qquad
    \widehat\mu=(L^t)_\#\widehat\nu.
\]
Thus $\Lambda$ records the intersections of the orbit $\varphi(\R^d)$ with $M$, while its spectrum is obtained by projecting the Fourier support of $\nu$ under $L^t$. This is the intersection--projection mechanism:
\[
\begin{tikzcd}[column sep=large, row sep=large]
    \begin{gathered}
        \mu\\[-2pt]
        \Lambda\subset\R^d
    \end{gathered}
    \arrow[r, "\mathcal A_\varphi"]
    \arrow[d, "\mathcal F_{\R^d}"']
    &
    \begin{gathered}
        \nu\\[-2pt]
        M\subset\T^n
    \end{gathered}
    \arrow[d, "\mathcal F_{\T^n}"]
    \\
    \begin{gathered}
        \widehat\mu\\[-2pt]
        S(\Lambda)\subset\Gamma\subset\R^d
    \end{gathered}
    &
    \begin{gathered}
        \widehat\nu\\[-2pt]
        \supp(\widehat\nu)\subset\Z^n
    \end{gathered}
    \arrow[l, "(L^t)_\#"']
\end{tikzcd}
\]
Here $\mathcal A_\varphi$ denotes the averaged pushforward
\[
    \mathcal A_\varphi(\mu)
    =
    \lim_{R\to\infty}
    \frac{\varphi_\#(\mu|_{B_R})}
         {\operatorname{vol}(B_R)}
    =\nu,
\]
with weak convergence of measures on $\T^n$.}
\subsubsection{Dimension of the Zariski closure of $\varphi(\Lambda)$}
In Section \ref{sec:trigsvanish}, using the same interpolation result as above and an ergodicity argument, we study the growth rate of the dimension of the space of Laurent polynomials that vanish on $M$ and for which all exponent vectors are contained in a growing bounded region of $\R^n$. A Hilbert polynomial type of argument then lets us deduce that the (complex) dimension of the Zariski closure $X:=\overline{M}^{\mathrm{Zar}}$ is equal to the (real) dimension of $M$ which is $c=n-d$. 
\subsubsection{Holonomicity leads to lighthouse}
As a consequence, we can express the measure $\nu$ on $M$ as integration of a rational differential top-form on a semi-algebraic subset of the smooth points of $X$. Classical results by Tougeron \cite{tougeron} then imply that the measure $\nu$ is holonomic, considered as a distribution on $\T^n$. In Section \ref{sec:holonomic} we prove that holonomicity together with the condition that the image of $\supp(\hat\nu)$ under $L^t\colon\Z^n\to\R^d$ is locally finite, implies that $\supp(\hat\nu)$ actually avoids an open cone that contains $\ker(L^t)\setminus\{0\}$. We use the so-called $b$-function of a holonomic ideal \cite{bfunctions} which gives many finite recurrence relations on the Fourier coefficients of $\nu$. This implies that $M$ is a lighthouse manifold, a notion we introduce and study in Section~\ref{sec:lighthouse}. Proposition~\ref{prop:lighthousegivesFQ} then supplies the previously unassumed polynomial growth of the Fourier transform. 
\subsubsection{Algebraicity of lighthouse manifolds}
In Section~\ref{sec:lighthouseimpliesalgebraic} we prove that every lighthouse manifold is of the form $M=X\cap\T^n$ where $X$ is a closed algebraic subvariety of $(\C^*)^n$. 
\subsubsection{Lee--Yang realization}
Finally, in Section~\ref{sec:final} we verify that $X$ satisfies an amoeba condition from \cite{alon2025higher}*{Theorem 8.4} that lets us deduce that $\Lambda$ comes from the Lee--Yang construction in \cite{alon2025higher}.

\section{{Uniformly discrete crystalline measures are of toral type}}\label{sec:toraltype}
It is inherent to the constructions in \cites{alon2025higher,lawtontsikh} that the spectrum of each resulting Fourier quasicrystal is contained in a finitely generated subgroup of $\R^d$. Following \cite{lawton-inverse}, we call this property \emph{toral type}. In this section we prove that the same conclusion holds {in fact for every crystalline measure with uniformly discrete support. In the following, we let $\mu$ be a crystalline measure on $\R^d$, let $\Lambda$ be its support and assume that $\Lambda$ is $\delta$-uniformly discrete for $\delta>0$.}

The first result of this section is 
\begin{thm}\label{thm: FQ toral}
    There exists a constant $C_d>0$ depending only on $d$ such that, for any 
    {crystalline measure $\mu$ on $\R^d$ with $\delta$-uniformly discrete support $\Lambda$ and Fourier support $S=\supp(\hat\mu)$, we have that $\mu$ is Bohr almost periodic and
    \[
        S\subset
        \spann_\Z\bigl(S\cap B_{R}\bigr),\qquad R=C_d/\delta.
    \]
    In particular, $\mu$ is of toral type.}
\end{thm}
{
\begin{ex}
    The assumption of being uniformly discrete in Theorem \ref{thm: FQ toral} cannot be removed. Indeed, there are crystalline measures $\mu$ (with unit weights) whose support is not even contained in a finitely generated $\Q$-vector space \cite{kurasovsarnak}*{Lemma 1}. Then the conclusion of Theorem \ref{thm: FQ toral} does not hold for the crystalline measure $\hat\mu$.
\end{ex}}
We need the following statement about measures supported on uniformly discrete sets, which generalizes \cite{lev2017fourier}*{Corollary 3.4} and might be of independent interest.
\begin{thm}\label{thm:measuresequal}
    There exists a constant $C=C_d$ depending only on the dimension $d$ such that the following holds.
    Let $\delta>0$ and $\mu_1,\mu_2$ be tempered (complex) measures on $\R^d$ such that $\supp(\mu_1)$ and $\supp(\mu_2)$ are both $\delta$-uniformly discrete. If $\hat{\mu}_1$ and $\hat{\mu}_2$ agree on a ball $B_R$ of radius $R={C/\delta}$, then $\mu_1=\mu_2$.
\end{thm}
Note that for $\mu_2=0$ this is \cite{lev2017fourier}*{Corollary 3.4}.
\subsection{Proof of Theorem \ref{thm:measuresequal}}
We will crucially rely on the following interpolation result that essentially goes back to Ingham \cite{ingham} in the one-dimensional case, and Kahane \cite{kahane} in the multi-dimensional case. 
{
\begin{lem}\label{lem:interpolation}
    There is a constant $K=K_d>0$ which depends on the dimension $d$ only, such that for every $\delta$-uniformly discrete set $\Lambda\subset\R^d$ and every $\lambda\in\Lambda$ there exists a Schwartz function $\varphi=\varphi_{\Lambda,\lambda}\in \mathcal{S}(\R^d)$ such that
    \[\varphi(x)=\begin{cases}
        1 & x=\lambda\\
        0 & x\in\Lambda\setminus\{\lambda\}
    \end{cases},\quad\text{and}\quad \supp(\widehat{\varphi})\subset B_{r},\]
    with $r={K}/{\delta}$. Furthermore, we can choose the $\varphi_{\Lambda,\lambda}$ such that there is a constant $M=M_\Lambda>0$ with $\|\varphi_{\Lambda,\lambda}\|_1\leq M$ for all $\lambda\in\Lambda$.
\end{lem}
\begin{proof}
    Let $\Lambda\subset\R^d$ be $\delta$-uniformly discrete for $\delta>0$.
    For a compact set $\Omega\subset\R^d$ denote by $\cB(\Omega)$ the Bernstein space of all bounded and continuous functions $f$ on $\R^d$ for which the distribution $\hat{f}$ is supported on $\Omega$. By \cite{lev2017fourier}*{Theorem 3.3} there is a constant $c>0$ which depends only on the dimension $d$ such that for every $\lambda\in\Lambda$ there exists $f_\lambda\in\cB(B_{c/\delta})$ with $f_\lambda(\lambda)=1$ and $f_\lambda(\lambda')=0$ for all $\lambda'\in\Lambda\smallsetminus\{\lambda\}$. As pointed out at the end of \cite{olevskii2012multi}*{Section 1.2} there further exists a constant $m>0$ that depends only on $\Lambda$, but not on $\lambda$, such that $\|f_\lambda\|_\infty\leq m$.\footnote{This can be seen by applying the open mapping theorem to the map $\cB(B_{c/\delta})\to\ell^\infty(\Lambda)$ given by restriction to $\Lambda$, which is surjective by \cite{lev2017fourier}*{Theorem 3.3}.}

    Now we can proceed as in the proof of \cite{lev2017fourier}*{Corollary 3.4}. Namely, for $\lambda\in\Lambda$ let $\psi_\lambda$ be a Schwartz function with $\psi_\lambda(\lambda)=1$ such that the support of $\hat\psi_\lambda$ is contained in $B_{1/\delta}$. Note that we can choose the $\psi_\lambda$ such that there exists $N>0$ with $\|\psi_\lambda\|_1\leq N$ for all $\lambda\in\Lambda$. Then $\varphi:=f_\lambda\psi_\lambda$ is a Schwartz function satisfying the desired conditions with $M=N\cdot m$. Indeed, compact Fourier support implies that $f_\lambda$ is $C^\infty$ with bounded derivatives, see for example \cite{littlewoodpaley}*{Lemma 2.1}. Thus the product with a Schwartz function is Schwartz.
\end{proof}

\begin{cor}\label{cor:transbounded}
    Let $\mu$ be a crystalline measure with uniformly discrete support. Then $\mu$ is \emph{translation-bounded} in the sense that
    \begin{equation*}
        \sup_{x\in\R^d}\int_{B_1(x)}|d\mu|<\infty.
    \end{equation*}
\end{cor}
\begin{proof}
    Let $\delta>0$ such that the support $\Lambda$ of $\mu$ is $\delta$-uniformly discrete and write $\mu=\sum_{\lambda\in\Lambda}a_\lambda\delta_\lambda$. By uniform discreteness, it suffices to find a uniform upper bound on the $|a_\lambda|$. Further write $\hat\mu=\sum_{\xi\in\supp(\hat\mu)}b_\xi\delta_\xi$. Let $K,M>0$ and $\varphi_{\Lambda,\lambda}$ be as in Lemma \ref{lem:interpolation}. Then
    \begin{equation*}
        |a_\lambda|=|\sum_{x\in\Lambda}a_x\varphi_{\Lambda,\lambda}(x)|=|\sum_{\xi\in\supp(\hat\mu)}b_\xi\hat\varphi_{\Lambda,\lambda}(-\xi)|\leq\sum_{\xi\in\supp(\hat\mu)\cap B_{K/\delta}}|b_\xi|\cdot M.
    \end{equation*}
    The right-hand side does not depend on $\lambda$ and, as $\supp(\hat\mu)$ is locally finite,  is finite. This proves the claim.
\end{proof}

\begin{cor}\label{cor:crystallineBAP}
    Let $\mu$ be a crystalline measure with uniformly discrete support. Then $\mu$ is a Bohr almost periodic measure.
\end{cor}
\begin{proof}
    This follows from Corollary \ref{cor:transbounded} and \cite{meyer2018global}*{Lemma 5.8}.
\end{proof}
}

We first use Lemma \ref{lem:interpolation} to prove that the supports of the measures in Theorem \ref{thm:measuresequal} must be close in the Hausdorff distance. We denote the Hausdorff distance between closed sets $E, F\subset\R^d$ by $d_{H}(E,F)$.
\begin{lem}\label{lem:closeinhausdorff}
    Let $\mu_1,\mu_2$ be tempered measures on $\R^d$ such that $\Lambda_1:=\supp(\mu_1)$ and $\Lambda_2:=\supp(\mu_2)$ are both $\delta$-uniformly discrete. 
    Let $r=K/\delta$ be the radius from Lemma \ref{lem:interpolation}.
    If  $\hat{\mu}_1$ and $\hat{\mu}_2$ agree on $B_{2r}$, then $d_{H}(\Lambda_1,\Lambda_2)\le \delta$.
\end{lem}
\begin{proof}
    Assume for the sake of a contradiction that the statement is false. Without loss of generality, we can assume that there exists $\lambda_2\in\Lambda_2$ whose distance to $\Lambda_1$ is larger than $\delta$. Then $\Lambda'=\Lambda_1\cup\{\lambda_2\}$ is $\delta$-uniformly discrete. Applying Lemma \ref{lem:interpolation} to $\Lambda'$ gives a Schwartz function $\varphi_1$ with $\varphi_1(\lambda_2)=1$, $\varphi_1(\lambda)=0$ for $\lambda\in\Lambda_1$, and $\supp(\widehat{\varphi}_1)\subset B_{r}$. Applying Lemma \ref{lem:interpolation} to $\Lambda_2$ gives another Schwartz function $\varphi_2$ with $\varphi_2(\lambda_2)=1$, $\varphi_2(\lambda)=0$ for $\lambda\in\Lambda_2\smallsetminus\{\lambda_2\}$, and $\supp(\widehat{\varphi}_2)\subset B_{r}$. Then $f=\varphi_1\cdot\varphi_2$ is a Schwartz function with $f(\lambda_2)=1$, $f(\lambda)=0$ for $\lambda\in(\Lambda_1\cup\Lambda_2)\smallsetminus\{\lambda_2\}$, and $\supp(\widehat{f})\subset B_{2r}$. Thus
    \begin{equation*}
        0=\int \overline f\,d\mu_1=\widehat\mu_1(\overline{\widehat f})=\widehat\mu_2(\overline{\widehat f})=\int\overline f\,d\mu_2=\mu_2(\{\lambda_2\})\neq0,
    \end{equation*}
    which is a contradiction.
\end{proof}
\begin{proof}[Proof of Theorem \ref{thm:measuresequal}]
    Let $K$ be the constant from Lemma \ref{lem:interpolation} and let $C=4K$. Let $\delta>0$ and let $r=K/\delta$ and $R=C/\delta=4r$. Let $\mu_1,\mu_2$ be tempered measures on $\R^d$ such that $\Lambda_1=\supp(\mu_1)$ and $\Lambda_2=\supp(\mu_2)$ are $\delta$-uniformly discrete and such that $\hat{\mu}_1$ and $\hat{\mu}_2$ agree on $B_R$. We have to show that $\mu_1=\mu_2$. To this end, let $\lambda_1\in\Lambda_1$. 
    Applying Lemma \ref{lem:closeinhausdorff} with separation parameter $\delta/2$ gives $d_{H}(\Lambda_1,\Lambda_2)\le \delta/2$, so there exists $\lambda_2\in\Lambda_2$ with $\|\lambda_1-\lambda_2\|\leq\frac{\delta}{2}$. 
    Because $\Lambda_1$ is $\delta$-uniformly discrete, the set $\Lambda_3=(\Lambda_1\smallsetminus\{\lambda_1\})\cup\{\lambda_2\}$ is $\frac{\delta}{2}$-uniformly discrete. We also note that $\lambda_1\not\in\Lambda_2$ and $\lambda_2\not\in\Lambda_1$, unless $\lambda_1=\lambda_2$. 
    Applying Lemma \ref{lem:interpolation} to $\Lambda_2$ gives a Schwartz function $\varphi_1$ with $\varphi_1(\lambda_2)=1$, $\varphi_1(\lambda)=0$ for $\lambda\in\Lambda_2\smallsetminus\{\lambda_2\}$, and $\supp(\widehat{\varphi}_1)\subset B_{r}$.   
    Applying Lemma \ref{lem:interpolation} to $\Lambda_3$ gives a Schwartz function $\varphi_2$ with $\varphi_2(\lambda_2)=1$, $\varphi_2(\lambda)=0$ for $\lambda\in\Lambda_1\smallsetminus\{\lambda_1\}$, and $\supp(\widehat{\varphi}_2)\subset B_{2r}$.
    Then $\rho=\varphi_1\cdot\varphi_2$ is Schwartz with $\rho(\lambda_2)=1$, $\rho(\lambda)=0$ for all $\lambda\in(\Lambda_1\cup\Lambda_2)\smallsetminus\{\lambda_1,\lambda_2\}$, and $\supp(\hat\rho)\subset B_{3r}$. Now
    \begin{equation*}
        \mu_1(\{\lambda_1\})\cdot\bar\rho(\lambda_1)=\int\bar{\rho}(x)d\mu_1(x)=\hat{\mu}_1\left(\bar{\hat{\rho}}\right)=\hat{\mu}_2\left(\bar{\hat{\rho}}\right)=\int\bar{\rho}(x)d\mu_2(x)=\mu_2(\{\lambda_2\}),
    \end{equation*}
    which shows that $\bar{\rho}(\lambda_1)=\frac{\mu_2(\{\lambda_2\})}{\mu_1(\{\lambda_1\})}$. 
    We note that, if $\lambda_1=\lambda_2$, this implies $\mu_1(\{\lambda_1\})=\mu_2(\{\lambda_1\})$ because $\rho(\lambda_2)=1$. If $\lambda_1\neq\lambda_2$, then there is a Schwartz function $g$ with $\supp(\hat{g})\subset B_{r}$ and $g(\lambda_1)\neq g(\lambda_2)$. Then $f=g\cdot\rho$ is Schwartz with $f(\lambda)=0$ for all $\lambda\in(\Lambda_1\cup\Lambda_2)\smallsetminus\{\lambda_1,\lambda_2\}$ and $\supp(\hat{f})\subset B_{R}$. We have
    \begin{align*}
        \mu_1(\{\lambda_1\})\bar f(\lambda_1)&=\int\bar{f}(x)d\mu_1(x)=\hat{\mu}_1\left(\bar{\hat{f}}\right)=\hat{\mu}_2\left(\bar{\hat{f}}\right)=\int\bar{f}(x)d\mu_2(x)\\&=\mu_2(\{\lambda_2\})\bar f(\lambda_2).
    \end{align*}
    By our choice of $f$, this implies
    \begin{equation*}
        \bar g(\lambda_2)=\bar f(\lambda_2)=\bar f(\lambda_1)\cdot\frac{\mu_1(\{\lambda_1\})}{\mu_2(\{\lambda_2\})}=\bar g(\lambda_1)\cdot\bar\rho(\lambda_1)\cdot\frac{\mu_1(\{\lambda_1\})}{\mu_2(\{\lambda_2\})}=\bar g(\lambda_1)
    \end{equation*}
    which contradicts our choice of $g$, and thus the case $\lambda_1\neq\lambda_2$ cannot occur. We have therefore shown that for every $\lambda_1\in\Lambda_1$ we have $\mu_1(\{\lambda_1\})=\mu_2(\{\lambda_1\})$. Applying the same argument with $\mu_1$ and $\mu_2$ interchanged proves $\mu_1=\mu_2$.
\end{proof}
\subsection{Proof of Theorem \ref{thm: FQ toral}}
We first state the following lemma.
{\begin{lem} \label{lem: AP and compactness}
  Let $\mu$ be a crystalline measure on $\R^d$ with $\delta$-uniformly discrete support $\Lambda$ and
  let $(\tau_j)_{j\in\N}$ be a sequence in $\R^d$.
  There is a measure $\mu'$ that is crystalline and whose support $\Lambda'$ is $\delta$-uniformly discrete, such that
  $\supp(\hat\mu')=\supp(\hat\mu)$, 
  and there exists a subsequence $(\tau_{j_k})_{k\in\N}$ so that for any $\xi\in\supp(\hat\mu)$
  \[\widehat{\mu'}(\{\xi\})=\lim_{k\to\infty}e^{-2\pi i \langle \tau_{j_k},\xi\rangle}\widehat{\mu}(\{\xi\}).\]
\end{lem}
\begin{proof}
By Corollary \ref{cor:crystallineBAP} the measure $\mu$ is Bohr almost periodic.
We denote by $\mu_j$ the translation of $\mu$ by $\tau_j$. We further let $\Omega(\mu)$ be the \emph{translational hull} of $\mu$, i.e., the closure of the set of all translations of $\mu$ by elements of $\R^d$, taken in the vague topology. By \cite{leelenz}*{Lemma 2} it is a compact set. Thus, after passing to a subsequence of $(\tau_j)_{j\in\N}$, the $\mu_j$ converge in the vague topology to some $\mu'$. By Lemma \ref{lem:vagueSchwartz} the $\mu_j$ and the $\hat\mu_j$ converge to $\mu'$ and $\hat\mu'$ as tempered distributions. The statement on the Fourier transform follows. Furthermore, for every $\epsilon>0$ and every $\lambda\in\supp(\mu')$, we have $\supp(\mu_j)\cap B_\epsilon(\lambda)\neq\emptyset$ for large enough $j$. Since each $\supp(\mu_j)$ is $\delta$-uniformly discrete, the set $\supp(\mu')$ is also $\delta$-uniformly discrete.
\end{proof}}

With this lemma, we can prove the next lemma.
{
\begin{lem}\label{lem: translations}
    Let $\mu$ be a crystalline measure with $\delta$-uniformly discrete support and fix arbitrary $R>0$. If there exists $\xi\in\supp(\hat\mu)$ such that $\xi\notin \spann_{\Z}(\supp(\hat\mu)\cap B_R)$, then there is a crystalline measure $\mu'$ with $\delta$-uniformly discrete support,  such that
    \[
        \hat\mu|_{B_R}
        =\hat\mu'|_{B_R},
        \qquad
        \hat\mu(\{\xi\})
        \ne\hat\mu'(\{\xi\}).
    \]
    In particular, $\mu\ne\mu'$.
\end{lem}}
\begin{proof}
    Let {$\{\omega_1,\ldots,\omega_N\}=\supp(\hat\mu)\cap B_R$} and let $\Gamma_0=\spann_\Z\{\omega_1,\ldots,\omega_N\}$. We first claim that there exists a sequence $\tau_j\in\R^d$ such that
    \[
        \lim_{j\to\infty}e^{2\pi i\langle\tau_j,\omega_i\rangle}=1
        \quad\text{for all $i, \quad$  and}
        \quad
        \lim_{j\to\infty}e^{2\pi i\langle\tau_j,\xi\rangle}= a_\xi\ne1.
    \]
This claim is a standard argument, based on Kronecker's approximation theorem. To prove it, let $H\subset (\R/\Z)^{N+1}$ be the closure of the image of the homomorphism
    \[
        T\colon\R^d\longrightarrow(\R/\Z)^{N+1},\qquad
        \tau\longmapsto
        (\langle\tau,\omega_1\rangle,\ldots,
        \langle\tau,\omega_N\rangle,\langle\tau,\xi\rangle).
    \]
    Denote the set of integer linear relations among $\omega_1,\ldots,\omega_N,\xi$ by 
    \[
        \mathcal R=\left\{(m,q)\in\Z^N\times\Z\mid
        \sum_{i=1}^Nm_i\omega_i+q\xi=0\right\}.
    \]
    By Kronecker's approximation theorem, $H$, the torus closure of the image of $T$, has the form 
    \[
        H=\{y\in(\R/\Z)^{N+1}\mid
        \langle(m,q),y\rangle=0\pmod{\Z}
        \text{ for all }(m,q)\in\mathcal R\}.
    \]
    Let $G_{\xi}$ be the subgroup of $\Z$ defined by
    \[G_{\xi}:=\{q\in\Z\ \mid\ \exists (m,q)\in\mathcal{R}\}=\{q\in\Z\ \mid\ q\xi\in\spann_{\Z}(\omega_{1},\ldots,\omega_N)\},\]
    so that for any  $\alpha\in\R/\Z$,
    \[(0,\ldots,0,\alpha)\in H\iff q\alpha=0\mod\Z,\ \text{ for all } q\in G_{\xi}.\]
 If $G_{\xi}=\{0\}$, we may choose $\alpha=1/2$. Otherwise $G_{\xi}=q_0\Z$ for some $q_0\in\N$, and the assumption that $\xi\notin\spann_{\Z}(\omega_1,\ldots,\omega_N)$ implies $q_0\ge2$. In such a case, take $\alpha=1/q_0$ as it satisfies $\alpha q=0\mod\Z$ for any $q\in G_{\xi}$. In both cases $(0,\ldots,0,\alpha)\in H$ and $a_\xi=e^{2\pi i\alpha}\ne1$. This proves the claim.

{Now the lemma follows by applying Lemma \ref{lem: AP and compactness} to the sequence $(\tau_j)_{j\in\N}$ from the claim.}
\end{proof}
\begin{proof}[Proof of Theorem \ref{thm: FQ toral}]
{The crystalline measure $\mu$ is Bohr almost periodic by Corollary \ref{cor:crystallineBAP}. Now, 
    choose $C=C_d$ as in Theorem \ref{thm:measuresequal} and put $R=C/\delta$. If there exists $\xi\in\supp(\hat\mu)$ such that $\xi\notin \spann_{\Z}(\supp(\hat\mu)\cap B_R)$, then Lemma \ref{lem: translations} gives a different crystalline measure $\mu'$ with $\delta$-uniformly discrete support, such that 
    \[\hat\mu|_{B_R}
        =\hat\mu'|_{B_R}.\]
    This contradicts Theorem \ref{thm:measuresequal}. Since $\supp(\hat\mu)$ is locally finite, $\supp(\hat\mu)\cap B_R$ is finite, so  $\spann_{\Z}(\supp(\hat\mu))$ is finitely generated. }
\end{proof}

\subsection{Preliminaries on submanifolds of a torus}\label{sec:geoemtric_def}
\begin{Def}
    Let $V\in\Gr(d,n)$ and let $W=V^\perp\in\Gr(c,n)$, where $c=n-d$. A closed topological submanifold $M\subset(\R/\Z)^n$ of dimension $c$ is \emph{transverse to $V$} if it is locally the graph of a continuous function $W\to V$. Namely, for every $x\in M$ there are an open neighborhood $U\subset W$ of $0$ and a continuous function $\psi\colon U\to V$, with $\psi(0)=0$, such that
    \[
        \{x+w+\psi(w)\mid w\in U\}
    \]
    is an open neighborhood of $x$ in $M$. If $L\in\R^{n\times d}$ has rank $d$, we say that $M$ is transverse to $L$ if it is transverse to its column span. We use the same terminology on $\T^n$ through the identification $x\mapsto\exp(2\pi i x)$.
\end{Def}
A submanifold transverse to $L$ is naturally equipped with a measure.
\begin{Def}
    Let $M\subset\T^n$ be a closed topological submanifold transverse to a rank $d$ matrix $L\in\R^{n\times d}$. Its \emph{$L$-directed measure} is defined, for Borel subsets $A\subset M$, by
    \[
        \nu_L(A):=\lim_{\epsilon\to0}\epsilon^{-d}
        \mathrm{Haar}\bigl(\{z\exp(2\pi iLx)\mid z\in A,
        \ x\in[0,\epsilon]^d\}\bigr),
    \]
    where the Haar measure on $\T^n$ is normalized to have mass one.
\end{Def}
\begin{rem}\label{rem:dependenceonL}
    Up to a positive scalar factor, $\nu_L$ depends only on the column span $V$ of $L$. Thus the normalized measure depends only on $V$ and is called the \emph{$V$-directed measure}.
\end{rem}
\begin{lem}\label{lem : product structure}
    Let $M\subset\T^n$ be a closed topological submanifold transverse to $L$. There is $\epsilon_0>0$ such that
    \[
        H\colon M\times[-\epsilon_0,\epsilon_0]^d\longrightarrow\T^n,
        \qquad H(z,x)=z\exp(2\pi iLx),
    \]
    is a homeomorphism onto a closed neighborhood of $M$. Its restriction to $M\times(-\epsilon_0,\epsilon_0)^d$ has open image. In these coordinates the Haar measure is $d\nu_L(z)\,dx$.
\end{lem}
\begin{proof}
    Put $V=\operatorname{im}(L)$ and choose $P\in\R^{n\times c}$ whose columns form a basis of $V^\perp$. In a local graph chart we may write
    \[
        z(u)=z_0\exp\bigl(2\pi i(Pu+L\psi(u))\bigr).
    \]
    Then $H(z(u),x)$ is represented in linear coordinates by $(u,\psi(u)+x)$. This is a homeomorphism onto an open set when the parameters range over an open set. Moreover, by Fubini's theorem the pullback of the Haar measure is $|\det[P\ L]|\,du\,dx$: for each fixed $u$, the change in the second coordinate is a translation. This proves the existence and local expression of $\nu_L$ in the preceding definition.

    Compactness gives a uniform $\epsilon_0$ for injectivity. Indeed, otherwise there would be distinct pairs $(z_j,x_j)$ and $(z'_j,x'_j)$ with $x_j,x'_j\to0$ and $H(z_j,x_j)=H(z'_j,x'_j)$. Passing to a subsequence, both $z_j$ and $z'_j$ converge to the same point of $M$, contradicting injectivity in a local graph chart. The map on the closed product is a homeomorphism onto its image by compactness, and the local description proves the remaining assertions.
\end{proof}
\subsection{Lawton's manifold}\label{sec:lawton}
By Theorem \ref{thm: FQ toral} every $\Lambda$ satisfying Assumption \ref{ass:crystalline} is a uniformly discrete Bohr almost periodic set of toral type. In \cite{lawton-inverse} Lawton associated to every such set a closed submanifold of the torus from which the set can be recovered. 
We only state this result for such sets with locally finite Fourier support, which is less general than in \cite{lawton-inverse}, and use the topological formulation of transversality given above. A direct proof is given in Appendix~\ref{app:lawton}.
\begin{thm}[Lawton's manifold]\label{thm: Lawton}
    Let $\Lambda\subset \R^d$ be a set satisfying Assumption \ref{ass:crystalline}, with Fourier support $S(\Lambda)\subset\R^{d}$ and Fourier coefficients $(c_{\xi})_{\xi\in S(\Lambda)}$.
 Let $\Gamma=\spann_{\Z}(S(\Lambda))$, let $n=\rank(\Gamma)$ (which is finite by Theorem \ref{thm: FQ toral}), and let $L\in\R^{n\times d}$ be a matrix such that $L^{t}:\Z^n\to\Gamma$ is an isomorphism. Define the homomorphism $\varphi\colon\R^d\to\T^n$ by $\varphi(x):=\exp(2\pi i Lx)$. Then, the Euclidean closure of the image of $\Lambda$ under $\varphi$ 
\begin{equation*}
    M:=\overline{\varphi(\Lambda)}\subset\T^n
\end{equation*}
is a submanifold of $\T^n$ with $\dim(M)=c=n-d$ that is transverse to $L$, such that:
\begin{enumerate}
    \item $\Lambda=\Lambda(M,L):=\varphi^{-1}(M)$.
    \item The averages of the push-forward of $\mu_{\Lambda}|B_{R}$ converge, as $R\to\infty$, to the $L$-directed measure $\nu:=\nu_{L}$ in the weak topology 
    \[
    \mathcal A_\varphi(\mu)
    :=
    \lim_{R\to\infty}
    \frac{\varphi_\#(\mu|_{B_R})}
         {\operatorname{vol}(B_R)}
    =\nu,
\]
    \item The Fourier transform of $\nu$ is given by 
\begin{equation*}
    \hat\nu(k)=\begin{cases}
        0,&\textrm{if }L^tk\not\in S(\Lambda), \\
        c_{L^tk},&\textrm{if }L^tk\in S(\Lambda),
    \end{cases} 
\end{equation*}
for all $k\in\Z^n$.
\end{enumerate}
\end{thm}
In a diagram, Theorem \ref{thm: Lawton} states that 
\[
\begin{tikzcd}[column sep=large, row sep=large]
    \begin{gathered}
        \mu\\[-2pt]
        \Lambda\subset\R^d
    \end{gathered}
    \arrow[r, "\mathcal A_\varphi"]
    \arrow[d, "\mathcal F_{\R^d}"']
    &
    \begin{gathered}
        \nu\\[-2pt]
        M\subset\T^n
    \end{gathered}
    \arrow[d, "\mathcal F_{\T^n}"]
    \\
    \begin{gathered}
        \widehat\mu\\[-2pt]
        S(\Lambda)\subset\Gamma\subset\R^d
    \end{gathered}
    &
    \begin{gathered}
        \widehat\nu\\[-2pt]
        \supp(\widehat\nu)\subset\Z^n
    \end{gathered}
    \arrow[l, "(L^t)_\#"']
\end{tikzcd}
\]

\begin{rem}
 Note that the choice of $\Z$-basis of $\Gamma$ identifies $\T^n$ with the Pontryagin dual $\Hom(\Gamma,\T)$ of $\Gamma$.  Under this identification, the submanifold $M$ of $\T^n$ corresponds to the closure of
 \begin{equation*}
    \{\varphi_x\in\Hom(\Gamma,\T)\mid x\in\Lambda\}
\end{equation*}   
 where $\varphi_x$ is the group homomorphism
\begin{equation*}
\varphi_x\colon\Gamma\to\T,\,\xi\mapsto e^{2\pi i\langle x,\xi\rangle}.    
\end{equation*}
This is a basis independent, hence canonical, description of $M$.
\end{rem}
\section{Trigonometric polynomials vanishing on crystalline supports}\label{sec:trigsvanish}
{Throughout the section we consider a non-zero crystalline measure $\mu$ with $\delta$-uniformly discrete support $\Lambda$ and denote its spectrum by $S=\supp(\hat\mu)$.
By Corollary \ref{cor:crystallineBAP} the measure $\mu$ is Bohr almost periodic. Note that this implies that $\Lambda$ is a Delone set. Furthermore, by Theorem \ref{thm: FQ toral} there is a finitely generated subgroup $\Gamma$ of $\R^d$ which contains $S$.}
Theorem \ref{thm:measuresequal} applied to $\mu_1=\mu$ and $\mu_2=0$ gives:
\begin{lem}\label{lem:constantC}
    There exists $C>0$ such that for a trigonometric polynomial $p$,
    \[p\ \text{  vanishes on $\Lambda$ }\ \iff \ \widehat{p\cdot\mu}\ \text{ vanishes on $B_C(0)$}.\]
\end{lem}
Let $L_1,\ldots,L_n$ be a $\Z$-basis of $\Gamma$ such that $L_1,\ldots,L_d$ are linearly independent over $\R$ and $L\in\R^{n\times d}$ be the matrix whose rows are $L_1,\ldots,L_n$.
Our goal is the following:
\begin{thm}\label{thm:dimzarclosoflawton}
    { Let $X\subset(\C^*)^n$ be the Zariski closure of the set
    \begin{equation*}
        \{\exp(2\pi iLx)\mid x\in\Lambda\}.
    \end{equation*}
    The (complex) dimension of {every irreducible component of} $X$ is at most $c=n-d$.
    In particular, if 
    $\Lambda\subset\R^d$ satisfies Assumption~\ref{ass:crystalline} and $M\subset\T^n\subset(\C^*)^n$ is Lawton's manifold, as defined in Theorem \ref{thm: Lawton}, then the (complex) dimension of {every irreducible component of} the Zariski closure $X$ of $M$ is equal to the (real) dimension of $M$ which is $c=n-d$.}
\end{thm}

Since $L_{1},\ldots,L_{d}$ are linearly independent over $\R$, up to change of coordinates, we can and will assume that $L_1,\ldots,L_d$ are the standard basis of $\R^d$, i.e., that
\begin{equation*}
L=\begin{pmatrix}
    \mathrm{id_{d}}\\ T
\end{pmatrix},\quad T\in\R^{c\times d}  .
\end{equation*}

\subsection{Equidistribution of $\Gamma/\Z^d$}
Let $\pi\colon \R^d\to\R^d/\Z^d$ be the quotient map, and recall that $\Gamma=L^{t}\Z^n=\Z^d+T^t\Z^c$, so $\Z^d$ is a subgroup of $\Gamma$. We define the subgroup $G\subseteq\R^d/\Z^d$
as the Euclidean closure of $\Gamma/\Z^d$ in the torus $\R^d/\Z^d$,
\begin{equation}
G:=\overline{\pi(\Gamma)}=\overline{\pi(T^t\Z^c)}.
\end{equation}
The goal in this section is to show that $\pi (\Gamma)$ equidistributes in $G$ according to the normalized uniform measure on $G$. 
We denote the normalized uniform measure on $G$ by $\mu_{G}$. We also define the subsets of $\Gamma$
\begin{equation*}
    \Gamma_{N}:=\{q+T^tm\mid q\in\Z^d,\ m\in\Z^c,\ \|m\|_{\infty}\le N\},\qquad \forall N\in\N,
\end{equation*}
and the associated counting measures of $\pi(\Gamma_{N})$
\begin{equation*}
    \mu_N:=\frac{1}{(2N+1)^{c}}\sum_{x\in \pi (\Gamma_{N})}\delta_{x}=\frac{1}{(2N+1)^{c}}\sum_{m\in \Z^c\cap[-N,N]^c}\delta_{\pi(T^tm)}
\end{equation*}
The equidistribution is 
\begin{equation}\label{eq:ergodiclimit}
    \int_Gfd\mu_G=\lim_{N\to\infty}\int_Gfd\mu_N
\end{equation}
for every bounded measurable function $f\colon G\to\R$ whose discontinuity set has zero $\mu_G$ measure. This is a standard ergodicity statement, see for example \cite{mullerrichard}*{Theorem 2.16} together with the Portmanteau theorem.

We apply \eqref{eq:ergodiclimit}  to a special class of functions. For any bounded Borel set $U\subset \R^d$ let $\mathbf{1}_U$ be its indicator function and define $I_U\colon G\to\R$ by 
\begin{equation*}
    I_U(x)=\sum_{y\in\pi^{-1}(x)}\mathbf{1}_U(y).
\end{equation*} 
It is defined such that
\begin{equation*}
\int_GI_Ud\mu_N=\frac{1}{(2N+1)^{c}}\sum_{x\in \pi(\Gamma_{N})}I_{U}(x)=\frac{|\Gamma_{N}\cap U|}{(2N+1)^{c}}.
\end{equation*}
The function $I_U$ is locally constant off $\pi(\partial U)$. Since $U$ is bounded, the number of points in each fiber $\pi^{-1}(x)\cap\partial U$ is uniformly bounded. We may conclude that $\int_GI_{\partial U}d\mu_G=0\iff \mu_{G}(\pi(\partial U))=0$, and if this holds then the discontinuity set of $I_{U}$ has zero $\mu_{G}$ measure. We deduce
\begin{lem}\label{lem:ANestimate}
    If $U\subset\R^d$ is a bounded Borel set with $\mu_{G}(\pi(\partial U))=0$,  then
    \begin{equation*}
        \lim_{N\to\infty}\frac{|\Gamma_N\cap U|}{(2N+1)^{c}}=\int_GI_Ud\mu_G.
    \end{equation*}
\end{lem}
\begin{rem}\label{rem:goodballs}
For each fixed center $x\in\R^d$, all but countably many radii $R>0$ satisfy $\mu_G(\pi(\partial B(x,R)))=0$. Indeed, the locally finite measure $U\mapsto\int_G I_U\,d\mu_G$ assigns positive mass to at most countably many of the disjoint spheres centered at $x$. Thus we may choose arbitrarily large admissible radii, and arbitrarily small admissible balls centered at zero. The latter have positive measure because their images contain a neighborhood of zero in $G$.
\end{rem}
\subsection{The ideal in the Laurent polynomial ring}\label{sec:idealinlaurent}
In the ring $\C[z_1^{\pm},\ldots,z_n^{\pm}]$ of Laurent polynomials we consider the ideal
\begin{equation*}
    I:=\{P\in \C[z_1^{\pm},\ldots,z_n^{\pm}]\mid P(\exp(2\pi i Lx))=0,\ \forall x\in\Lambda\},
\end{equation*}
of all Laurent polynomials $P$ such that the trigonometric polynomial
\begin{equation*}
    P(\exp(2\pi i Lx))=P\left(e^{2\pi i\langle L_1,x\rangle},\ldots,e^{2\pi i\langle L_n,x\rangle}\right)
\end{equation*}
vanishes on $\Lambda$. Let $X\subset(\C^*)^n$ be the common zero set of all polynomials in $I$. 

{In the case that $\mu$ has unit weights, meaning that $\Lambda$ satisfies Assumption \ref{ass:crystalline}}, if $M\subset \T^n$ is the manifold from Theorem \ref{thm: Lawton}, then 
\[M\subset X\cap \T^n.\]
This means that
\begin{equation*}
    \dim_{\C}(X)\ge\dim_{\R}(X\cap\T^n)\ge\dim_{\R}(M)=n-d,
\end{equation*}
since $\T^n$ is totally real in $(\C^*)^n$.
The goal is to prove:
\begin{thm}\label{thm:codimd}
    The (complex) dimension of {every irreducible component of}  $X$ is upper bounded by $ c=n-d$.
\end{thm}
For a nonzero Laurent polynomial $0\neq P\in\C[z_1^{\pm},\ldots,z_n^{\pm}]$ we use the notation $\newt(L,P):= L^t\newt(P)\subset\R^d$ where $\newt(P)\subset\R^n$ is the Newton polytope of $P$.
\begin{lem}\label{lem:boundeddiameter}
    There exists $r>0$ such that the diameter of $\newt(L,P)$ is at least $r$ for every $0\neq P\in I$.
\end{lem}
\begin{proof}
    Let $p(x)=P(\exp(2\pi i Lx))$ and let $R$ be such that $\Lambda$ is $R$-relatively dense. Namely, $\R^d=\cup_{x\in\Lambda}B(x,R)$. Let $D$ be the diameter of $\newt(L,P)$. We will prove that
    \begin{equation}\label{eq:bernsteinconclusion}
        \|p\|_\infty\leq\pi R D\|p\|_\infty
    \end{equation}
    which implies the claim for $r=\frac{1}{\pi R}$. Let $x\in\R^d\smallsetminus\Lambda$ and $y\in\Lambda$ a point with $\|x-y\|< R$. Define
    \begin{equation*}
        g(t):=p\left(y+t\cdot\frac{x-y}{\|x-y\|}\right),
    \end{equation*}
    so that $g(0)=0$ and $g(\|x-y\|)=p(x)$. By choosing suitable $a\in\R$, the function  $f(t)=e^{2\pi i at}g(t)=\sum a_j e^{i \omega_j t} $ is a finite sum of exponentials with frequencies $\omega_j\in[-\pi D,\pi D]$ for all $j$. Bernstein's inequality \cite{bernstein} implies that
    \begin{equation*}
        \|f'\|_\infty\leq\pi D\|f\|_\infty.
    \end{equation*}
    Since $f(\|x-y\|)=\int_{0}^{\|x-y\|}f'(t)dt$, and $\|x-y\|<R$, we conclude
    \begin{equation*}
        \|p(x)\|=\|f(\|x-y\|)\|\leq\int_0^R\|f'(t)\|dt\leq R \pi D\|f\|_\infty\leq R\pi D\|p\|_\infty,
    \end{equation*}
    which shows (\ref{eq:bernsteinconclusion}).
\end{proof}
\begin{lem}\label{lem:indepconst}
    Let $c=n-d$. Given any bounded subset  $U\subset\R^d$,  define
    \begin{equation*}
        V(U,N)=\spann_\C\{z^m\mid m\in\Z^n,\, L^tm\in U\cap\Gamma_N\},\qquad\forall N\in\N.
    \end{equation*}
    Then, there is a constant $C>0$, independent of the choice of $U$, such that
    \begin{equation*}
    \limsup_{N\to\infty}\frac{\dim(V(U,N)/(I\cap V(U,N)))}{(2N+1)^c}\leq C.
    \end{equation*}
\end{lem}
\begin{proof}
    Enlarge the constant $C$ from Lemma \ref{lem:constantC}, if necessary, using Remark~\ref{rem:goodballs}, so that the ball $B=B(0,C)\subset\R^d$ satisfies the condition $\mu_{G}(\pi(\partial B))=0$ in Lemma \ref{lem:ANestimate}. We consider the linear map
    \begin{equation*}
        \Psi\colon V(U,N)\to\cS'(\R^d),\quad \Psi(Q):=\widehat{q\mu}\cdot\mathrm{1}_{B},
    \end{equation*}
    where $q(x)=Q(\exp(2\pi iLx))$ and $\mathrm{1}_{B}$ is the indicator of $B$. By Lemma \ref{lem:constantC},
    \[\ker(\Psi)=I\cap V(U,N).\]
    We now want to bound the dimension of the image of $\Psi$. Since $\supp(\widehat{\mu})=S$ and  $\supp(\widehat{q})\subset U\cap\Gamma_{N}$ for any $Q\in V(U,N)$, one has
$$\supp(\Psi(Q))=\supp(\widehat{q\cdot\mu})\cap B\subset \left(S+(U\cap\Gamma_{N})\right)\cap B$$
Since $\Gamma_{N}\cap U$ is finite, and $S$ is locally finite, the set $\left(S+(U\cap\Gamma_{N})\right)\cap B$ is a finite subset of $\Gamma$. 
Notice that $x\in \left(S+(U\cap\Gamma_{N})\right)\cap B$ if and only if $x\in B$ and can be written as $x=s+\gamma$ for some $s\in S$ and $\gamma\in \Gamma_N\cap U$, so \[s=x-\gamma\in S\cap (B-U).\]
The set $B-U$ is bounded, so $S\cap (B-U)$ is a finite subset of $\Gamma$, and as such it is contained in $\Gamma_{N_U}$ for some constant $N_{U}\in\N$ that depends on $U$ but is independent of $N$. So $s\in \Gamma_{N_U}$ and $\gamma\in \Gamma_{N}$, so $x=s+\gamma\in\Gamma_N+\Gamma_{N_U}\subset\Gamma_{N+N_U}$. But we also know that $x\in B$, so 
\[\supp(\Psi(Q))\subset \Gamma_{N+N_{U}}\cap B\qquad \forall Q\in V(U,N).\]
We conclude that 
\begin{equation*}
        {\dim(V(U,N)/(I\cap V(U,N)))}=\dim(\mathrm{image}(\Psi))\leq |\Gamma_{N+N_U}\cap B|,
    \end{equation*}
for all $N$, with $N_{U}$ fixed. 
    
    Now Lemma \ref{lem:ANestimate} implies 
    \[\limsup_{N\to\infty}\frac{{\dim(V(U,N)/(I\cap V(U,N)))}}{(2N+1)^c}\le \lim_{N\to\infty}\frac{|\Gamma_{N+N_U}\cap B|}{(2N+1)^c}=\int_{G}I_{B}d\mu_{G},\]
    and the limit $C:=\int_{G}I_{B}d\mu_{G}$ is independent of the choice of $U$.
\end{proof}
\begin{lem}\label{lem:sumset}
    Let $X\subset(\C^*)^n$ be a closed irreducible algebraic subvariety of dimension $e\geq1$ and $\C[X]$ be its coordinate ring, i.e., the ring of functions $X\to\C$ given by the restriction of a Laurent polynomial to $X$.
    There is a finite dimensional $\C$-vector space $W_0\subset \C[X]$ such that every finite dimensional $\C$-vector space $W_0\subset W\subset \C[X]$ satisfies
    \begin{equation*}
        \dim(W^k)\geq\dim(W)\cdot\binom{e+k-1}{k-1}-(k-1)\binom{e+k-1}{k}
    \end{equation*}
     for all $k\geq1$. Here we denote by $W^k$ the span of all products of $k$ elements of $W$.
\end{lem}
\begin{proof}
    Let $x\in X$ be a smooth point of $X$ and $a_1,\ldots,a_e\in\C[X]$ elements which vanish at $x$ such that $da_1(x),\ldots,da_e(x)$ are linearly independent. By the analytic implicit function theorem we can, locally around $x$, write every element of $\C[X]$ as a power series in $a_1,\ldots,a_e$. Hence we get a map of $\C$-algebras $\C[X]\to\C[[x_1,\ldots,x_e]]$ which sends $a_i$ to $x_i$. By the identity principle and because $X$ is irreducible, this map is injective. Set $W_0=\operatorname{span}_{\C}\{1,a_1,\ldots,a_e\}$. In the following, we identify $\C[X]$ with its image under this map. 
    Consider the map 
    \begin{equation*}
    v\colon\C[[x_1,\ldots,x_e]]\smallsetminus\{0\}\to\Z^e    
    \end{equation*}
    which sends a nonzero power series to the smallest occurring exponent with respect to the lexicographic monomial ordering. If $W\subset \C[X]$ is a finitely generated $\C$-vector space, then
    \begin{equation}\label{eq:prodtosum}
        \dim(W)=|v(W\smallsetminus\{0\})| \text{ and } k\cdot v(W\smallsetminus\{0\})\subset v(W^k\smallsetminus\{0\}),
    \end{equation}
    where $k\cdot K$ for $K\subset\R^e$ denotes the set of sums of $k$ elements of $K$. Since $v(W_0\smallsetminus\{0\})$ contains an affine basis of $\R^e$, we can apply \cite{sumsets}*{Theorem 1.5} and obtain for all finite dimensional $W_0\subset W\subset \C[X]$ that
    \begin{equation*}
        |k\cdot v(W\smallsetminus\{0\})|\geq |v(W\smallsetminus\{0\})|\cdot\binom{e+k-1}{k-1}-(k-1)\binom{e+k-1}{k},
    \end{equation*}
    which implies the claim by (\ref{eq:prodtosum}).
\end{proof}
\begin{lem}
    If 
    \begin{enumerate}
    \item $I$ is a prime ideal,
        \item there is a bounded subset $U_0\subset\R^d$  with $\int_G I_{U_0}\,d\mu_G>0$ and $\mu_G(\pi(\partial U_0))=0$, that has the property that $V(U_0,N)\cap I=\{0\}$ for every $N\geq 1$ ($V(U_0,N)$ is defined as in Lemma \ref{lem:indepconst}), and
        \item there exists a constant $C>0$ such that for every bounded $U\subset\R^d$ we have
    \begin{equation*}
        \limsup_{N\to\infty}\frac{\dim((V(U,N)+I)/I)}{(2N+1)^c}\leq C,
    \end{equation*}
    \end{enumerate}
    then the (complex) dimension of $X$ is at most $c$.
\end{lem}
\begin{proof}
    Let $A=\C[z_1^{\pm},\ldots,z_n^{\pm}]/I$ be the coordinate ring of $X$ and let $e\geq1$ be its Krull dimension. Further, let $W_0\subset A$ be as in Lemma \ref{lem:sumset} and for $N\geq1$ define
    \begin{equation*}
        W_N=(V(U_0,N)+I)/I+W_0\subset A.
    \end{equation*}
    Since $V(U_0,N)\cap I=\{0\}$, Lemma \ref{lem:ANestimate} implies that
    \begin{equation}\label{eq:limitWN}
        \lim_{N\to\infty}\frac{\dim W_N}{(2N+1)^{c}}=\int_GI_{U_0}d\mu_G.
    \end{equation}
    Now let $N_1\in\N$ and $U_0\subset U_1\subset\R^d$ bounded such that $W^{'}_{N_1}:=(V(U_1,N_1)+I)/I$ contains $W_0$. With these choices and for $N\geq N_1$ we have 
    \begin{equation*}
        W_N\subset (V(U_1,N)+I)/I.
    \end{equation*}
    Thus for $k\in\N$ we have that
    \begin{equation*}
        W_N^k\subset (V(kU_1,kN)+I)/I,
    \end{equation*}
    where $kU_1$ is the $k$-fold Minkowski sum of $U_1$.
    Lemma \ref{lem:sumset} then implies
    \begin{equation*}
        \limsup_{N\to\infty}\frac{1}{(2N+1)^c}\left(\dim(W_N)\cdot\binom{e+k-1}{k-1}-(k-1)\binom{e+k-1}{k}\right)\leq k^c C.
    \end{equation*}
    Then, using (\ref{eq:limitWN}), for all $k\in\N$ we get
    \begin{equation*}
        \binom{e+k-1}{k-1}\cdot\int_GI_{U_0}d\mu_G\leq k^c C.
    \end{equation*}
    Since the left-hand side is a polynomial in $k$ of degree $e$, the claim $e\leq c$ follows.
\end{proof}
\begin{cor}
    If $I$ is prime, then $X$ has dimension at most $c$.
\end{cor}
\begin{proof}
    We apply the previous lemma.
    By Remark \ref{rem:goodballs} there are arbitrarily small balls $B$ satisfying the conditions $\int_GI_{B}d\mu_G>0$ and $\int_GI_{\partial B}d\mu_G=0$. By Lemma \ref{lem:boundeddiameter}, choosing the ball $B$ to be small enough also ensures the condition $V(B,N)\cap I=\{0\}$ for every $N\geq 1$ from the previous lemma. Finally, the constant $C$ required by the previous lemma exists by Lemma \ref{lem:indepconst}.
\end{proof}
\begin{proof}[Proof of Theorem \ref{thm:codimd}]
    Let $\fp_1,\ldots,\fp_l$ be the minimal prime ideals containing $I$. Their zero sets $X_1,\ldots,X_l$ are the irreducible components of $X$ and it suffices to show that each $X_i$ has codimension at least $d$. To this end, choose $F_j\in\fp_j\smallsetminus\fp_i$ for $j\neq i$.  Then the product $G=\prod_{j\neq i}F_j$ vanishes identically on each $X_j$, $j\neq i$, but not on $X_i$. Let $g=G(\exp(2\pi i Lx))$ be the corresponding trigonometric polynomial on $\R^d$. The measure $\mu'=g\cdot\mu$ and its support $\Lambda'=\supp(\mu')$ still satisfy the conditions from the beginning of the section. 
    Moreover, the corresponding ideal of Laurent polynomials is exactly $\fp_i$, so that the previous corollary applies.
\end{proof}
\begin{proof}[Proof of Theorem \ref{thm:dimzarclosoflawton}]
    It follows from Theorem \ref{thm:codimd} that the dimension of {every irreducible component of}  $X$ is at most $c$. {In the case $\mu$ has unit weights, since} $X$ is the Zariski closure of the $c$-dimensional submanifold $M$ of $\T^n$, the other inequality follows; {see also the discussion at the beginning of Section \ref{sec:idealinlaurent}.}
\end{proof}
\subsection{Consequences for Lawton's manifold}\label{sec:lawton2}
Consider the situation of Theorem \ref{thm:dimzarclosoflawton}. Namely, $\Lambda\subset\R^d$ satisfies Assumption~\ref{ass:crystalline}, $M\subset\T^n$ is Lawton's manifold, as defined in Section \ref{sec:lawton}, and $X\subset(\C^*)^n$ is its Zariski closure, which has {pure} dimension $c=n-d$ by Theorem \ref{thm:dimzarclosoflawton}. Let $X_{\rm sm}$ be the smooth part of $X$. Then $M_{\rm sm}:=M\cap X_{\rm sm}$ is a real-analytic set. Furthermore, it is a union of connected components of $X_{\rm sm}(\T):=X_{\rm sm}\cap\T^n$: it is relatively closed and by invariance of domain also relatively open. 

This shows in particular that $M_{\rm sm}$ is a semialgebraic subset of $\T^n$ and $M\smallsetminus M_{\rm sm}$ is contained in a semialgebraic subset of dimension less than $c$. 
From Lawton's description of the measure $\nu=\nu_L$ \cite{lawton-inverse} it follows that  one can orient $M_{\rm sm}$ in a way so that $\nu$ is given by integrating on $M_{\rm sm}$ against a regular differential form $\omega_L$ that we now define. Given $I=\{i_1,\ldots,i_c\}$ with $0<i_1<\cdots<i_c\leq n$ we consider the $c$-form
\begin{equation*}
    \omega_I=\frac{1}{(2\pi i)^c}\frac{dz_{i_1}\wedge\cdots\wedge dz_{i_c}}{z_{i_1}\cdots z_{i_c}}
\end{equation*}
and the sign $ s(I)=(-1)^{\sum_{i\in I}i-\sum_{i=1}^ci} $. Then
\begin{equation}\label{eq:topform}
    \omega_L:=\sum_{I\in\binom{[n]}{c}}L_{I}\cdot s(I)\cdot\omega_I,
\end{equation}
where $L_I$ is the determinant of the $d\times d$ submatrix of $L$ obtained by deleting the rows indexed by $I$ and $\binom{[n]}{c}$ is the set of $c$-element subsets of $[n]:=\{1,\ldots,n\}$.
\section{Lighthouse distributions and manifolds}\label{sec:lighthouse}
Let $n=c+d$ and denote by $\Gr(c,n)$ the Grassmannian of real $c$-dimensional subspaces of $\R^n$.
\begin{Def}
    Let $\nu$ be a distribution on $\T^n$ and $V\in\Gr(c,n)$. We say that $\nu$ \emph{strongly avoids $V$} if there exists an open neighborhood $U\subset\Gr(c,n)$ of $V$ such that for all $V'\in U$ we have: 
    \begin{equation*}
        V'\cap\supp(\widehat{\nu})\subset\{0\}.
    \end{equation*}
\end{Def}

\begin{Def}\label{def:lighthousemanifold}
    Let $M$ be a closed submanifold of $\T^n$ of dimension $c$, let $V\in\Gr(c,n)$ and $W\in\Gr(d,n)$ be the orthogonal complement of $V$. We say that $M$ is a \emph{lighthouse manifold with respect to $W$} if:
    \begin{enumerate}
        \item $M$ is transverse to $W$, and
        \item the $W$-directed measure $\nu_W$ on $M$ strongly avoids $V$.
    \end{enumerate}
    If $L\in\R^{n\times d}$ is a matrix of rank $d$, then we say that $M$ is a \emph{lighthouse manifold with respect to $L$} if it is one with respect to the column span of $L$.
    Finally, we say that a closed submanifold $M\subset\T^n$ is a \emph{lighthouse manifold} if it is a lighthouse with respect to some linear subspace of dimension $d=n-\dim(M)$.
\end{Def}
\begin{ex}\label{ex:leeyanglighthouse}
    Let $X\subset(\C^*)^n$ be a strict Lee--Yang variety of dimension $c$. Then its torus part $X(\T)=X\cap\T^n$ is a lighthouse manifold with respect to any $V=L\R^d$ with $L\in\R^{n\times d}$ that has only positive $d\times d$ minors, by \cite{alon2025higher}*{Lemma 7.4 and Theorem 3.3}.
\end{ex}
The following lemmas will be useful later on.
\begin{lem}\label{lem:lighthouselemma1}
    Let $S\subset\R^n$ and $V\in\Gr(c,n)$ be the left-kernel of a matrix $L\in\R^{n\times d}$. The following are equivalent:
    \begin{enumerate}
        \item There exists $C>0$ such that for all $k\in S$ we have $\|k\|\leq C\cdot\|L^tk\|$.
        \item There exists an open neighborhood $U\subset\Gr(c,n)$ of $V$ such that for all $W\in U$ we have $W\cap S\subset\{0\}$.
    \end{enumerate}
\end{lem}
\begin{proof}
    We first prove the direction $(1)\Rightarrow(2)$. For any linear subspace $W\subset\R^n$ we denote by $P_W\in\R^{n\times n}$ the matrix representing the orthogonal projection onto $W$. Consider the set
    \begin{equation*}
        U=\{W\in\Gr(c,n)\mid C\|L^tP_W\|_{\rm op}<1\}
    \end{equation*}
    where $\|\cdot\|_{\rm op}$ denotes the operator norm. We will prove that $U$ satisfies the conditions in $(2)$. The set is clearly open and, since $L^tP_V=0$, it contains $V$. Now let $W\in U$ and $k\in W\cap S$. We have
    \begin{equation*}
        \|k\|\leq C\cdot\|L^tk\|=C\cdot\|L^tP_Wk\|\leq \underbrace{C\cdot\|L^tP_W\|_{\rm op}}_{<1}\cdot\|k\|,
    \end{equation*}
    which implies $k=0$. 
    For proving the direction $(2)\Rightarrow(1)$ we assume that $(1)$ is not true. This means that there is a sequence $(k_i)_{i\in\N}\subset S\smallsetminus\{0\}$ such that $\lim_{i\to\infty}\frac{\|L^tk_i\|}{\|k_i\|}=0$. We can assume that, after passing to a subsequence, the $\frac{k_i}{\|k_i\|}$ converge to a vector $v_1$ on the unit sphere. This $v_1$ is necessarily in the kernel of $L^t$ and thus in $V$. We complete it to a basis $v_1,\ldots,v_c$ of $V$. Consider the vector spaces $W_i$ spanned by $k_i,v_2,\ldots,v_c$. For sufficiently large $i$ their dimension is $c$ and they converge in $\Gr(c,n)$ to $V$. Since each $W_i$ intersects $S\smallsetminus\{0\}$ in $k_i$, this shows that $(2)$ does not hold.
\end{proof}
\begin{lem}\label{lem:lighthouselemma2}
    Let $S\subset\Z^n$, $V\in\Gr(c,n)$ and $L\in\R^{n\times d}$ satisfy the conditions of Lemma \ref{lem:lighthouselemma1}. Let $P\in\R^{c\times n}$ whose row-span $V'$ is in the open neighborhood $U$ of $V$ from part $(2)$ of Lemma \ref{lem:lighthouselemma1}. Then the following holds:
    \begin{enumerate}
        \item The projected set $L^tS\subset\R^d$ is locally finite.
        \item For $m\in\Z^n$ there are at most finitely many $l\in\R^c$ such that $P^tl+m\in S$.
    \end{enumerate}
\end{lem}
\begin{proof}
    Part (1) follows immediately from part (1) of Lemma \ref{lem:lighthouselemma1} and the fact that $S\subset\Z^n$ is locally finite in $\R^n$. For proving part (2) let $Q\in\R^{n\times d}$ be a matrix whose left-kernel is $V'$. Since $U$ is an open neighborhood of $V'$, we can apply Lemma~\ref{lem:lighthouselemma1} to the matrix $Q$ and get a constant $C>0$ such that for all $k\in S$ we have $\|k\|\leq C\cdot \|Q^tk\|$. Now if $k=P^tl+m\in S$, then
    \begin{equation*}
        \|k\|\leq C\cdot \|Q^tk\|=C\cdot \|Q^t(P^tl+m)\|=C\cdot \|Q^tm\|.
    \end{equation*}
    Therefore, for fixed $m$ there are only finitely many possibilities for $k\in S\subset\Z^n$, and thus, since $P^t$ is injective, for $l$.
\end{proof}
Lighthouse manifolds can be used to produce Fourier quasicrystals.
\begin{prop}\label{prop:lighthousegivesFQ}
Let $L\in\R^{n\times d}$ have rank $d$, and let
$M\subset\T^n$ be a lighthouse manifold
with respect to $L$. Let $\nu$ be its $L$-directed measure, and denote 
\[
    \Lambda:=\Lambda(M,L),
    \qquad
    S=L^t\supp\widehat\nu.
\]
Then $S$ is discrete (i.e., locally finite), and  $\Lambda$ is a Delone FQ. In particular, $\Lambda$ satisfies Assumption \ref{ass:crystalline}, with  
\[
    \widehat{\mu_\Lambda}
    =
    \sum_{\xi\in S}c_\xi\delta_\xi,
    \qquad
    c_\xi=\sum_{k\in F_\xi}\widehat\nu(k),
\]
and the sum defining $c_{\xi}$ is over the fiber $ F_\xi=
    \{k\in\supp\widehat\nu\ \mid\ L^tk=\xi\}$ which is finite. 
In particular,
\[
    S(\Lambda)=\{\xi\in S\ \mid\ c_\xi\ne0\}.
\]
If, in addition, the rows of $L$ are linearly independent
over $\Q$, then
\[
    S(\Lambda)=S,\qquad
    c_{L^tk}=\widehat\nu(k),
\]
and
\[
    M=\overline{\{\exp(2\pi iLx):x\in\Lambda\}}.
\]
\end{prop}
\quash{\begin{prop}\label{prop:lighthousegivesFQ}
    Let $L\in\R^{n\times d}$ be a matrix of rank $d$ whose rows are linearly independent over $\Q$. Let $M\subset\T^n$ be a lighthouse manifold with respect to $L$ and let $\nu=\nu_L$ be its $L$-directed measure. Then $\Lambda=\Lambda(M,L)$ is a Delone Fourier quasicrystal whose spectrum is $S(\Lambda)=L^t\supp(\hat\mu)$ with Fourier coefficients $c_{L^tk}=\hat\mu(k)$ for $k\in\Z^n$. Finally, we have
    \begin{equation*}
        M=\overline{\{\exp(2\pi i Lx) \mid x\in\Lambda\}}\subset\T^n.
    \end{equation*}
\end{prop}}
\begin{proof} 
By Lemma \ref{lem : product structure}, we have an $\epsilon_0>0$ and a closed  neighborhood of $M$, call it $M_{0}\subset \T^n$, which is homeomorphic to $M\times [-\epsilon_0,\epsilon_0]^d$, with the homeomorphism $$H(z,x)=z\cdot\exp(2\pi i Lx),$$
and the Haar measure of $\T^n$ restricted to $M_0$ is $d\nu(z)dx$.
Let
$\rho\in C_{c}^{
\infty
}(\R^d)$ be non-negative, with $\supp(\rho)\subset Q_{d}:=[0,1]^d$ and $\hat{\rho}(0)=\int\rho dx=1$. For $\epsilon>0$, let $\rho_{\epsilon}(x):=\epsilon^{-d}\rho(\frac{x}{\epsilon})$, so that
\begin{equation*}
    \supp(\rho_{\epsilon})\subset \epsilon Q_{d},\quad \widehat{\rho_{\epsilon}}(\xi)=\hat{\rho}(\epsilon \xi),
\end{equation*} and in particular, as $\epsilon\to 0$, the following weak limit holds
\begin{equation*}
    \rho_{\epsilon}(x)dx\to\delta_{0},\qquad\text{as tempered distributions}.
\end{equation*}
For $\epsilon\in(0,\epsilon_0)$ define  
$$g_{\epsilon}\colon \T^n\to \R,\quad g_{\epsilon}(z)=0\ \forall z\notin M_{0},\quad g_{\epsilon}(H(z,x)):=
    \rho_{\epsilon}(x).$$
Notice that $g_{\epsilon}$ is continuous, and has Fourier coefficients
\begin{align*}
\widehat{g_{\epsilon}}(k) & :=\int_{\T^n}z^{-k}g_{\epsilon}(z)d\mathrm{Haar}(z)\\ &=\int_{z\in M,x\in[-\epsilon_{0},\epsilon_0]^d}z^{-k}e^{-2\pi i \langle L^{t}k,x\rangle}\rho_{\epsilon}(x)d\nu(z)dx\\ &=\widehat{\nu}(k)\widehat{\rho}(\epsilon L^{t}k).    
\end{align*}
In fact, $g_{\epsilon}\in C^{\infty}$. To see why, we use 
\begin{enumerate}
    \item $|\widehat{\nu}(k)|\le \widehat{\nu}(0)$ for all $k$, since $\nu$ is non-negative and finite,
    \item since $\rho\in C_{c}^{\infty}(\R^d)$, then for any $N\in\N$ there exists $C_N>0$ such that $$|\hat{\rho}(\epsilon L^t k)|\le C_{N}(1+\epsilon^2\|L^t k\|^2)^{-N},\quad \forall k\in\Z^n,$$
    \item $\|k\|\le C\|L^t k\|$ for all $k\in\supp(\widehat{\nu})$, by the lighthouse definition.
\end{enumerate}
We conclude that for any $N\in\N, \epsilon\in(0,\epsilon_0)$, there exists $C_{N,\epsilon}>0$ such that $|\widehat{g_{\epsilon}}(k)|\le C_{N,\epsilon}(1+\|k\|^2)^{-N}$ for all $k\in\Z^n$. We conclude that $g_{\epsilon}\in C^{\infty}(\T^n)$, and in particular its Fourier series converges absolutely and uniformly in $z\in\T^n$, 
$$g_{\epsilon}(z)=\sum_{k\in\supp(\widehat{\nu})}\widehat{\nu}(k)\widehat{\rho}(\epsilon L^{t}k)z^k.$$
Define the coefficients $(c_{\xi})_{\xi\in S}$, with $S=L^t\supp(\widehat{\nu})$, using the fibers  
\[F_{\xi}=\{k\in\supp\widehat{\nu}\ \mid\ L^tk=\xi\},\quad c_{\xi}:=\sum_{k\in F_{\xi}}\hat{\nu}(k), \quad \forall \xi\in S.\]
By Lemma~\ref{lem:lighthouselemma1}, there exists $C>0$
such that $\|k\|\le C\|\xi\|$ for every
$k\in F_{\xi}$, so there exists $C'>0$ such that, letting $S_R:=S\cap B_R$,
 \[\sum_{\xi\in S_{R}}|F_{\xi}|\le \left|\Z^n\cap B_{CR}\right|\le C'(1+R)^n,\quad \forall R>0.\]
 Using $|\widehat{\nu}(k)|\le \widehat{\nu}(0)$ for all $k$, and a dyadic decomposition argument, we get that 
 \begin{equation}\label{eq: fibersum}
 \sum_{\xi\in S}\frac{|c_{\xi}|}{(1+\|\xi\|)^m} \le \widehat{\nu}(0)\sum_{\xi\in S}\frac{|F_{\xi}|}{(1+\|\xi\|)^m}<\infty,\quad \forall m>n+1.    
 \end{equation}

Let $\mu=\mu_{\Lambda}$ be the counting measure of $\Lambda$. By construction, $\Lambda$ is $\epsilon_0$-uniformly discrete, so $\mu$ is tempered. Define  $$f_{\epsilon}\colon\R^d\to \R,\quad f_{\epsilon}(x):=g_{\epsilon}(\exp(2\pi i Lx)).$$
It inherits the Fourier series from $g_{\epsilon}$,
$$f_{\epsilon}(x)=\sum_{k\in\supp(\widehat{\nu})}\widehat{\nu}(k)\widehat{\rho}(\epsilon L^{t}k)e^{2\pi i \langle L^t k,x\rangle}=\sum_{\xi\in S }c_{\xi}\widehat{\rho}(\epsilon \xi)e^{2\pi i \langle \xi,x\rangle},$$
that converges absolutely and uniformly in $x\in\R^d$, since
\[\sum_{\xi\in S }|c_{\xi}\widehat{\rho}(\epsilon \xi)e^{2\pi i \langle \xi,x\rangle}|\le C_{N,\epsilon}\sum_{\xi\in S }\frac{|c_{\xi}|}{(1+\|\xi\|)^N}<\infty, \]
for $N>n+1$.
It follows that 
$f_{\epsilon}$ is Bohr almost periodic.
Therefore, for any Schwartz function $h\in\mathcal{S}(\R^d)$, on the one hand,

$$ I_{\epsilon}(h):=\int_{\R^d}f_{\epsilon}(x)\hat{h}(x)dx= \sum_{\xi\in S}c_{\xi}\widehat{\rho}(\epsilon \xi)h(\xi). $$
On the other hand, by definition, 
$$f_{\epsilon}(x)\hat{h}(x)=\sum_{y\in\Lambda}\rho_{\epsilon}(x-y)\hat{h}(x).$$
For each $y\in\Lambda$,
$\int_{\R^d}\rho_\epsilon(x-y)\widehat h(x)\,dx\to\widehat h(y)$.
Moreover, the Schwartz decay of $\widehat h$ and compact support of $\rho$ lead to a constant $C_h>0$, independent of $0<\epsilon<\epsilon_0$, such that
\[
    \int_{\R^d}|\rho_\epsilon(x-y)\widehat h(x)|\,dx
    \le C_h(1+\|y\|)^{-d-1},
\]
which is summable over the uniformly discrete
set $\Lambda$. Therefore, we can change the order of summation, integration, and limit, so that
$$\lim_{\epsilon\to 0} I_{\epsilon}(h)=\sum_{y\in\Lambda}\lim_{\epsilon\to 0}\int_{\R^d} \rho_{\epsilon}(x-y)\hat{h}(x)dx=\sum_{x\in\Lambda}\widehat{h}(x).$$
Finally, since $h$ is Schwartz, we have new constants $C_{N,h}$ such that $|h(\xi)|\le C_{N,h}(1+\|\xi\|)^{-N}$. Using $\sup_{\xi}|\widehat{\rho}(\xi)|=\widehat{\rho}(0)=1$, and taking $N>n+1$, there is a bound uniform in $\epsilon$, 
\[\sum_{\xi\in S}|c_{\xi}\widehat{\rho}(\epsilon \xi)h(\xi)|\le C_{N,h}\sum_{\xi\in S}\frac{|c_{\xi}|}{(1+\|\xi\|)^N}<\infty, \]
so we can replace the order of the limit and sum, 
\begin{align*}
    \int \widehat{h}d\mu=\sum_{x\in\Lambda}\widehat{h}(x)
    = & \lim_{\epsilon\to 0} I_{\epsilon}(h)\\
    = & \sum_{\xi\in S}\lim_{\epsilon\to 0}c_{\xi}\widehat{\rho}(\epsilon \xi)h(\xi)\\
    = & \sum_{\xi\in S}c_{\xi}h(\xi).
\end{align*}
We conclude that 
$$\widehat{\mu}=\sum_{\xi\in S}c_{\xi}\delta_{\xi},$$and it is strongly tempered, since $$|\widehat{\mu}|(B_{R})=\sum_{\xi\in S_{R}}|c_{\xi}|\le \widehat{\nu}(0) \sum_{\xi\in S_{R}}|F_{\xi}| \le C'(1+R)^n.$$
To see that $\Lambda\ne\emptyset$, namely that $\mu$ is not zero, we note that $\widehat{\mu}(\{0\})=\widehat{\nu}(\{0\})\ne0$. 
This means that $\Lambda$ is a uniformly discrete FQ, and as such it is also Delone.  

If in addition, the rows of $L$ are linearly independent over $\Q$ so $L^t:\Z^n\to \R^d$ is injective, $|F_{\xi}|=1$ for all $\xi\in S$, so in particular there is no cancellation and $c_{\xi}\ne 0 $ for all $\xi\in S$. Moreover, the orbit $\{\exp(2\pi i L x)\ \mid\ x\in\R^d\}$ is dense in $\T^n$. Since $M$ is a closed manifold transverse to $L$, then $M=\overline{\{\exp(2\pi i Lx) \mid x\in\Lambda\}}$.
\end{proof}
Note that the definition of a lighthouse manifold (Definition \ref{def:lighthousemanifold}) does not require any smoothness. However, smoothness is a consequence of the lighthouse condition.
\begin{prop}\label{prop:lighthousesmooth}
    Every lighthouse manifold $M\subset\T^n$ is $C^\infty$.
\end{prop}
\begin{proof}
    Let $c=\dim(M)$ and $d=n-c$. Choose $L\in\R^{n\times d}$ such that $M$ is a lighthouse with respect to $L$. Denote the translated manifold $M_{z}:=z^{-1}M$ for every $z\in\T^n$. Notice that $M_{z}$ is a lighthouse with respect to $L$ if and only if $M$ is, and it is $C^{\infty}$ if and only if $M$ is. Therefore, to prove the statement, it will be enough to assume that $\textbf{1}:=(1,\ldots,1)\in M$ and prove that $M$ is $C^{\infty}$ in a small neighborhood of $\textbf{1}$. We proceed with this assumption and goal.
    
    Let $\nu$ be the $L$ directed measure on $M$, and $\nu_{z}$ the $L$-directed measure on $M_{z}$, so $\hat\nu_{z}(k)=z^k\hat\nu(k)$ for all $k\in\Z^n$. Define $\Lambda=\Lambda(M,L)$, and for any $z\in\T^n$ define $\Lambda_z:=\Lambda(M_{z},L)$. By Proposition \ref{prop:lighthousegivesFQ} each $\Lambda_z$ is a Delone FQ, with spectrum contained in the set $S(\Lambda_z)\subset S=L^{t}\supp{\widehat{\nu}}$ independent of $z$,     
    and Fourier coefficients,
    \[\widehat{\mu_{\Lambda_z}}=\sum_{\xi\in S}c_{\xi}(z)\delta_{\xi},\qquad c_{\xi}(z)=\sum_{k\in F_{\xi}}z^k\widehat{\nu}(k),\]
    with the fibers $F_{\xi}:=\{k\in\supp{\widehat{\nu}}\mid  L^tk=\xi\}$. 

    Fix $O\subset \T^n$ a small neighborhood of $\textbf{1}$. By making $O$ small enough, there is a unique point $y(z)\in\Lambda_z$ which is closest to the origin among the points of $\Lambda_z$. We use $y(z)$ to define a parameterization of $M\cap O$. 
    
    Let $P\in\R^{n\times c}$ be a matrix whose columns are an orthonormal basis of $\ker(L^t)$, and define $z(x)=\exp(2\pi i Px)$ for $x\in\R^c$. Then, by taking $O$ to be smaller if needed, there is an open neighborhood of the origin $U\subset\R^c$, such that
    \begin{equation*}
        U\to\T^n,\, x\mapsto z(x)\cdot\exp(2\pi iLy(z(x)))
    \end{equation*}
    parametrizes $M\cap O$, the neighborhood of $\textbf{1}$.
    
For $i=1,\ldots,d$ we choose Schwartz functions $f_i\in\cS(\R^d)$ such that $\hat{f}_i$ are real and the support of $\hat{f}_i$ is a small neighborhood around the origin and such that on a smaller neighborhood it agrees with the $i$th coordinate function. Since these are finitely many Schwartz functions, there are constants $C_{N}>0$ such that $\max_{i}|f_{i}(\xi)|\le C_{N}(1+\|\xi\|)^{-N}$ for all $\xi\in\R^d$ and all $N\in\N$. 
Further, let
    \begin{equation*}
        \psi=\begin{pmatrix}
            \psi_1\\ \vdots \\
            \psi_d
        \end{pmatrix}\colon\T^n\to\R^d,\qquad \psi_i(z):=\sum_ {k\in\Z^n}z^k\hat\nu(k)f_i(L^tk)=\sum_{\xi\in S}c_{\xi}(z)f_{i}(\xi).
    \end{equation*}
The summation formula of $\Lambda_z$ gives
$\psi_i(z)=\sum_ {x\in\Lambda_{z}}\widehat{f_{i}}(x).$ 
Therefore, when $z$ is close enough to $\textbf{1}$, so that $\|y(z)\|$ is small enough, the defining features of $\widehat{f_i}$ give
\[\psi(z)=y(z).\]
Each $c_\xi$ is a finite Laurent polynomial. Positivity
and the lighthouse bound give, for every $\alpha\in\N_0^n$,
\[
    \sup_{z\in\T^n}|\partial^\alpha c_\xi(z)|
    \le C_\alpha|F_\xi|(1+\|\xi\|)^{|\alpha|}.
\]
Combining this with Schwartz decay and choosing
$N>n+|\alpha|+1$, we obtain
\[
    \sum_{\xi\in S}\sup_{z\in\T^n}
    |\partial^\alpha(c_\xi(z)f_i(\xi))|
    \le C_{\alpha,N}
    \sum_{\xi\in S}
    \frac{|F_\xi|}{(1+\|\xi\|)^{N-|\alpha|}}
    <\infty,
\]
by the fiber-summability estimate \eqref{eq: fibersum} in the proof
of Proposition~\ref{prop:lighthousegivesFQ}.
Thus all differentiated series converge absolutely and
uniformly. This proves that
$\psi_i\in C^\infty(\T^n)$ for all $i$, and therefore $y(z)$ is $C^\infty$ for $z$ in a small neighborhood of $\textbf{1}$.

 Since each $\psi_i$ is $C^\infty$, so is $M$.
\end{proof}
\section{Holonomic $D$-ideals}\label{sec:holonomic}
We recall some facts about holonomic $D$-ideals. The book \cite{dmodules} will serve as our main reference.
A \emph{$D$-ideal} is a left-ideal of the Weyl algebra $$D_n=\C\langle z_1,\ldots,z_n,\partial_1,\ldots,\partial_n\rangle.$$
Every element $p$ of $D_n$ can be uniquely written as
\begin{equation*}
    p=\sum_{(\alpha,\beta)\in E}c_{\alpha\beta}\cdot z^\alpha\partial^\beta
\end{equation*}
for a finite subset $E\subset\N_0^{2n}$ and non-zero $c_{\alpha\beta}\in\C$. Let $v\in\R^n$ and consider the following two types of \emph{initial forms} of such $p$:
\begin{align*}
    \ini_*(p)&:=\sum c_{\alpha\beta}\cdot z^\alpha w^\beta\in\C[z_1,\ldots,z_n,w_1,\ldots,w_n],\\
    \ini_v(p)&:=\sum c_{\alpha\beta}\cdot z^\alpha \partial^\beta\in\C\langle z_1,\ldots,z_n,\partial_1,\ldots,\partial_n\rangle,
\end{align*}
where the sums are over all pairs $(\alpha,\beta)\in E$ for which $\sum_{i=1}^n\beta_i$ (in the first sum) and $\langle v,\beta-\alpha\rangle$ (in the second sum) is maximal. Here $w_1,\ldots,w_n$ are additional commutative variables. The \emph{characteristic variety} $\textrm{ch}(I)\subset\C^{2n}$ of a $D$-ideal $I$ in $D_n$ is the common zero set of all $\ini_*(p)$ for $p\in I$.
\begin{Def}[\cite{dmodules}*{Definition 1.4.8}]
    A $D$-ideal $I\subset D_n$ is \emph{holonomic} if its characteristic variety has dimension $n$.
\end{Def}
The Weyl algebra $D_n$ acts on a distribution $T$ on $\T^n\subset(\C^*)^n$ by $(z_iT)(f):=T(z_if)$ and $(\partial_iT)(f):=-T(\partial_if)$.
We say that a distribution $\nu$ on $\T^n$ is \emph{holonomic} if its annihilator (a left-ideal of $D_n$) is holonomic. 
The following is the main result of this section.
\begin{thm}\label{thm:discreteimplieslighthouseforholonomic}
    Let $\nu$ be a holonomic distribution on $\T^n$ and let $S=\supp(\hat{\nu})\subset\Z^n$ be its Fourier support.
    Let $L\in\R^{n\times d}$ have $\Q$-linearly independent rows and full rank over $\R$. Let $V=\ker(L^t)\subset\R^n$ and assume that the set $L^tS\subset\R^d$ is locally finite. Then $\nu$ strongly avoids $V$.
\end{thm}

The proof of Theorem \ref{thm:discreteimplieslighthouseforholonomic} uses the so-called $b$-function of a holonomic $D$-ideal whose definition we now recall from \cite{dmodules}*{Chapter 5}. Let $v\in\R^n\setminus\{0\}$ and $I\subset D_n$ be a holonomic $D$-ideal. We define $\ini_v(I)$ to be the $\C$-span of all $\ini_v(p)$ for $p\in I$. It is again a $D$-ideal by \cite{dmodules}*{Corollary 1.1.2}. Denote $\theta_i=z_i\partial_i$ and consider the ideal
\begin{equation*}
    \ini_v(I)\cap\C[v_1\theta_1+\cdots+v_n\theta_n]
\end{equation*}
of the univariate polynomial ring $\C[s]$ where $s=\sum_{i=1}^nv_i\theta_i$. The generator $b(s)$ of this principal ideal is called the \emph{$b$-function of $I$ with respect to $v$}. Note that it is well-defined up to a non-zero multiplicative constant.
\begin{thm}[\cite{dmodules}*{Theorem 5.1.2}]
    Let $0\neq v\in\R^n$ and $I\subset D_n$ be a holonomic $D$-ideal. Then the $b$-function $b(s)$ of $I$ with respect to $v$ is not the zero polynomial.
\end{thm}

\begin{lem}\label{lem:recurrence}
    Let $0\neq v\in\R^n$ and $\nu$ be a holonomic distribution on $\T^n$. There are:
    \begin{enumerate}
        \item a finite set $F\subset\Z^n$ with $\langle m,v\rangle>0$ for all $m\in F$, and
        \item for each $m\in F$ a multivariate polynomial $q_m\in\C[z_1,\ldots,z_n]$,
    \end{enumerate}
    such that for all $k\in\Z^n$ the following recurrence relation holds:
    \begin{equation}\label{eq:rec}
        b(\langle k-\vec{1}, v\rangle)\cdot\hat{\nu}(k)+\sum_{m\in F}q_m(k-m-\vec{1})\cdot\hat{\nu}(k-m)=0,
    \end{equation}
    where $b$ is the $b$-function of the annihilator of $\nu$ with respect to $v$ and $\vec{1}\in\R^n$ is the all-ones vector.
\end{lem}

\begin{proof}
    Let $I$ be the annihilator of $\nu$. By definition of the $b$-function, there are $p_1,\ldots,p_N\in I$ such that
    \begin{equation*}
        b(\langle\theta,v\rangle)=\sum_{i=1}^N\ini_v(p_i).
    \end{equation*}
    Since on the left-hand side each monomial is of the form $z^\alpha\partial^\alpha$ for some $\alpha\in\N_0^n$, we can assume without loss of generality that every monomial in each $\ini_v(p_i)$ on the right-hand side is of the form $z^\alpha\partial^\beta$ with $\langle\beta-\alpha,v\rangle=0$. Then
    \begin{equation*}
        b(\langle\theta,v\rangle)=\sum_{i=1}^N\ini_v(p_i)=\ini_v\left(\sum_{i=1}^Np_i\right)=\ini_v(p)
    \end{equation*}
    where $p=\sum_{i=1}^Np_i\in I$. Then we can write 
    \begin{equation*}
        p=b(\langle\theta,v\rangle)+\sum_{(\alpha,\beta)\in E}c_{\alpha\beta}\cdot z^\alpha\partial^\beta=b(\langle\theta,v\rangle)+\sum_{(\alpha,\beta)\in E}c_{\alpha\beta}\cdot z^{\alpha-\beta}\prod_{i=1}^n\prod_{j=0}^{\beta_i-1}(\theta_i-j)
    \end{equation*}
    with $E\subset\N_0^{2n}$ some finite set with $\langle\alpha-\beta,v\rangle>0$ for all $(\alpha,\beta)\in E$. Hence the finite set $F=\{\alpha-\beta\mid(\alpha,\beta)\in E\}$ satisfies the first requirement of the claim. Now for $m\in F$ we define $$q_m=\sum_{(\alpha,\beta)\in E,\,\alpha-\beta=m}c_{\alpha\beta}\cdot \prod_{i=1}^n\prod_{j=0}^{\beta_i-1}(z_i-j)\in\C[z_1,\ldots,z_n].$$
    Then we have $p=b(\langle\theta,v\rangle)+\sum_{m\in F} z^{m}\cdot q_m(\theta).$
    To prove the recurrence relation (\ref{eq:rec}), we observe that 
    $\widehat{(\theta_i \nu)}(k)=(k_i-1)\hat{\nu}(k)$ 
     and $\widehat{z^m \nu}(k)=\hat{\nu}(k-m)$ for all $k,m\in\Z^n$. 
     Thus for $q\in\C[z_1,\ldots,z_n]$, we have
     \begin{equation*}
         \widehat{z^mq(\theta)\nu}(k)=q(k-m-\vec{1})\hat\nu(k-m)
     \end{equation*}
     Applying this to $p\cdot \nu=0$, we get the desired identity.
\end{proof}
We will combine the previous lemma with the following.
\begin{lem}\label{lem:descent}
    Let the rows of $L\in\R^{n\times d}$ be linearly independent over $\Q$. Let $S\subset\Z^n$ such that $L^tS\subset\R^d$ is locally finite. Assume that for all $v\in\R^d$ with $\|v\|=1$ there are $R_v>0$ and $F_v\subset\Z^n$ finite such that:
    \begin{enumerate}
        \item $m^tLv>0$ for all $m\in F_v$, and
        \item if $k\in S$ with $k^tLv>R_v$, then $k-m\in S$ for some $m\in F_v$.
    \end{enumerate}
    Then there exists $C>0$ and a finite set $T\subset\Z^n$ such that every $k\in S$ is a sum of at most $C\cdot\|L^tk\|$ many elements of $T$.
\end{lem}
\begin{proof}
    For $v\in\R^d$ with $\|v\|=1$ let $0<\delta_v<1$ such that $\frac{1}{2}m^tLv>\delta_v$ for all $m\in F_v$. Consider the following spherical neighborhood
    \begin{equation*}
        O_v=\{w\in\R^d\mid\|w\|=1,\,\langle w,v\rangle>\frac{1}{2},\,\forall m\in F_v\colon m^tLw>\delta_v\}
    \end{equation*}
    of $v$. By compactness of the sphere, there are finitely many $v_1,\ldots,v_r\in S^{d-1}$ such that the open sets $O_{v_1},\ldots,O_{v_r}$ cover $S^{d-1}$. Define
    \begin{equation*}
        F=F_{v_1}\cup\cdots\cup F_{v_r}.
    \end{equation*}
    We define:
    \begin{equation*}
           N=\max_{m\in F}\|L^tm\| \textrm{ and }
        \delta=\min_{i=1}^r\delta_{v_i},
    \end{equation*}
    and choose some $R\geq\max(2R_{v_1},\ldots,2R_{v_r}, N^2/\delta,1)$. If $F=\emptyset$, we define $N=0$. Furthermore, let
    \begin{equation*}
        B=\{m\in S\mid \|L^tm\|\leq R\} \textrm{ and }T=F\cup B.
    \end{equation*}
    Note that the set $T$ is finite because $L^tS$ is locally finite and since the homomorphism $\Z^n\to\R^d,\,m\mapsto L^tm$ is injective. Now let $k\in S$ with $\rho:=\|L^tk\|>R$ and $w:=L^tk/\rho$. There exists $i\in\{1,\ldots,r\}$ such that $w\in O_{v_i}$. Then we have
    \begin{equation*}
        k^tLv_i=\langle \rho w, v_i\rangle>\frac{\rho}{2}>\frac{R}{2}\geq R_{v_i}.
    \end{equation*}
    Therefore, by assumption (2), there exists $m\in F_{v_i}$ with $k':=k-m\in S$. We have
    \begin{align*}
        \|L^tk'\|^2&=\|\rho w-L^tm\|^2=\langle \rho w-L^tm,\rho w-L^tm\rangle \\
        &=\rho^2-\underbrace{2\rho(m^tLw)}_{>2\rho\delta}+\underbrace{\|L^tm\|^2}_{\leq N^2}\leq\rho^2-\rho\delta+\underbrace{(N^2-\rho\delta)}_{<0}\\
        &<\rho^2-\rho\delta<\rho^2-\rho\delta+\delta^2/4=(\rho-\delta/2)^2.
    \end{align*}
    This shows that $\|L^tk'\|<\|L^tk\|-\delta/2$. Iterating this, we arrive after at most 
    $$\max(0,\lceil2\cdot\frac{\|L^tk\|-R}{\delta}\rceil)$$
    steps of subtracting elements from $F$ at an element of $B$. This shows that every $k\in S\smallsetminus B$ is a sum of at most 
    \begin{equation*}
        2\cdot\frac{\|L^tk\|-R}{\delta}+2=\frac{2}{\delta}\cdot\|L^tk\|+\underbrace{2-\frac{2R}{\delta}}_{\leq0}\leq \frac{2}{\delta}\cdot\|L^tk\|
    \end{equation*}
    many elements of $T$. Since $B$ is finite, we get the claim by replacing $\frac{2}{\delta}$ by some large enough constant $C$.
\end{proof}
\begin{cor}\label{cor:lighthouseS}
    In the situation of Lemma \ref{lem:descent} there exists $\tilde{C}>0$ such that every $k\in S$ satisfies $\|k\|\leq\tilde{C}\cdot\|L^tk\|$.
\end{cor}
\begin{proof}
    This follows immediately from Lemma \ref{lem:descent} by choosing $\tilde{C}$ to be the product of $C$ and the maximal norm of an element of $T$.
\end{proof}
\begin{proof}[Proof of Theorem \ref{thm:discreteimplieslighthouseforholonomic}]
    By Corollary \ref{cor:lighthouseS} it suffices to verify the conditions of Lemma \ref{lem:descent}. Thus let $v\in\R^d$ with $\|v\|=1$ and define $F_v$ to be the finite set $F$ of Lemma \ref{lem:recurrence} applied to $0\neq Lv\in\R^n$. Property (1) of Lemma \ref{lem:recurrence} then implies condition (1) of Lemma \ref{lem:descent}. Furthermore, let $b_v$ be the $b$-function of the annihilator of $\nu$ with respect to $Lv$ and let $R_v$ be some positive number that is larger than every real root of $b_v(t-\langle\vec{1},Lv\rangle)$. Then the recurrence relation (\ref{eq:rec}) in Lemma \ref{lem:recurrence} implies condition (2) of Lemma \ref{lem:descent}.
\end{proof}

\subsection{Lawton's manifold is a lighthouse}
Now let $\Lambda\subset\R^d$ satisfy Assumption~\ref{ass:crystalline}, $L\in\R^{n\times d}$ a matrix whose rows form a basis of the group generated by the spectrum of $\Lambda$, $M\subset\T^n$ Lawton's manifold and $\nu=\nu_L$ the $L$-directed measure on it as defined in Section \ref{sec:lawton}.
As in Section \ref{sec:lawton2}, we can realize $\nu$ by integrating on the semialgebraic subset $M_{\rm sm}$ against the differential form $\omega_L$ from (\ref{eq:topform}). The following lets us apply Theorem \ref{thm:discreteimplieslighthouseforholonomic}.

\begin{prop}
    The distribution $f\mapsto\int_Mf\omega_L$ on $\T^n$ is holonomic.
\end{prop}

\begin{proof}
    Since $\T^n$ is compact, it is enough by \cite{tougeron}*{Proposition 1.5} to show that this distribution is locally holonomic in the sense that every point has an open neighborhood restricted to which the distribution is holonomic. Let $\pi\colon\tilde{X}\to X$ be a resolution of singularities together with an antiholomorphic involution $\tilde\sigma$ on $\tilde{X}$ so that $\pi\circ\tilde\sigma=\sigma\circ\pi$, where $\sigma$ is the coordinate-wise involution $z\mapsto1/\bar{z}$ on $X$ \cite{hironaka}. Let $\tilde{\omega}:=\pi^*\omega_L$ and $M_{\rm sm}=M\cap X_{\rm sm}$ where $X_{\rm sm}$ is the smooth part of $X$. We have
    \begin{equation*}
        \int_Mf\omega_L=\int_{M_{\rm sm}}f\omega_L=\int_{\pi^{-1}(M_{\rm sm})}(f\circ\pi)\tilde\omega.
    \end{equation*}
    Thus by \cite{tougeron}*{Proposition 2.13} it is enough to show that integrating against $\tilde\omega$ on the semialgebraic open subset $\pi^{-1}(M_{\rm sm})$ of the fixed points of $\tilde\sigma$ is locally holonomic. For this, by \cite{tougeron}*{Corollaire 2.10} it is enough to show that the density of $\tilde\omega$ is locally holonomic on $\pi^{-1}(M_{\rm sm})$ and this follows from \cite{tougeron}*{Corollaire 1.4}.
\end{proof}

\begin{cor}\label{cor:lawtonlighthouseS}
    The measure $\nu$ strongly avoids the kernel of $L^t$, i.e.~there exists $C>0$ such that every $k\in\supp(\hat\nu)$ satisfies $\|k\|\leq{C}\cdot\|L^tk\|$.
\end{cor}

\begin{proof}
    This follows from the previous proposition and Theorem \ref{thm:discreteimplieslighthouseforholonomic}.
\end{proof}

\begin{cor}\label{cor:lawtonlighthouseandsmooth}
    Lawton's manifold $M$ is a lighthouse with respect to $L$.
\end{cor}
\begin{proof}
The manifold is transverse to $L$ by Theorem~\ref{thm: Lawton}, and its directed measure strongly avoids $\ker(L^t)$ by Corollary~\ref{cor:lawtonlighthouseS}.
\end{proof}
\begin{cor}\label{cor:crystallineisFQ}
Every set satisfying Assumption~\ref{ass:crystalline} is a Delone Fourier quasicrystal.
\end{cor}
\begin{proof}
Let $M$ and $L$ be as in Theorem~\ref{thm: Lawton}. By Corollary~\ref{cor:lawtonlighthouseandsmooth}, $M$ is a lighthouse with respect to $L$. Proposition~\ref{prop:lighthousegivesFQ} applies because the rows of $L$ form a $\Z$-basis of the spectrum group. Since $\Lambda=\Lambda(M,L)$, it gives the claim, including the growth bound~\eqref{ass:growth}.
\end{proof}

\section{Lighthouse implies algebraic}\label{sec:lighthouseimpliesalgebraic}
\subsection{A general criterion for real varieties}
Let $f\colon X\to Y$ be a generically finite morphism of affine varieties over $\R$ whose restriction to every irreducible component of $X$ is dominant.
We will denote by $A=\C[Y]$ and $B=\C[X]$ the coordinate rings of (the complexifications of) $Y$ and $X$, respectively.
We further assume that $Y$ is smooth and irreducible. Let $M\subset X(\R)$ be a subset that is Zariski dense in $X$ such that the restriction $f|_M\colon M\to Y(\R)$ is a topological covering map. 
Let $K=\C(Y)$ be the complex function field of $Y$ and similarly $F_i=\C(X_i)$ where $X_1,\ldots,X_r$ are the irreducible components of $X$. Finally, we set $F=\C(X)=F_1\times\cdots\times F_r$.
The goal of this subsection is to prove the following.
\begin{thm}\label{thm:realfiberedcrit}
    Assume that there exists a nowhere-vanishing function $g_M\colon M\to\C$,
   such that for every $h\in\C[X]$ the function
    \begin{equation}\label{eq:defGh}
        Y(\R)\to\C,\, y\mapsto\sum_{x\in f^{-1}(y)\cap M}g_M(x)h(x)
    \end{equation}
    agrees with a regular function $G_h\in\C[Y]$. Then $f$ is finite, real-fibered, in the sense that $f^{-1}(Y(\R))=X(\R)$, and $M=X(\R)$. Moreover, every point of $X(\R)$ is a smooth point of $X$.
\end{thm}
Each $F_i$ is a finite and separable field extension of $K$. The ring $B=\C[X]$ is a finitely generated ring extension of $A$, torsion-free as $A$-module and $F=K\otimes_AB$.
Recall that the \emph{trace map} $\tr=\tr_{F/K}\colon F\to K$ sends $a\in F$ to the trace of
\begin{equation*}
    F\to F,\, x\mapsto a\cdot x,
\end{equation*}
and that the $K$-bilinear form
\begin{equation*}
    F\times F\to K,\, (x,y)\mapsto\tr_{F/K}(xyz)
\end{equation*}
is nondegenerate for all $z\in F^\times$ \cite{lang}*{VI Theorem 5.2}.
\begin{lem}\label{lem:tracealg}
    There exists $g\in\C(X)$ such that
    for every $h\in B$, we have $G_h=\tr(gh)$, and which agrees with $g_M$ on a dense Zariski open subset of $M$. Moreover, there is a Zariski dense open affine subset $U\subset Y$ such that $f^{-1}(U(\R))\subset M$.
\end{lem}
\begin{proof}
    Since $Y(\R)$ is Zariski dense in $Y$, the regular function $G_h$ is uniquely determined by (\ref{eq:defGh}). The map $B\to A,\, h\mapsto G_h$ is thus clearly $A$-linear. Since $F=K\otimes_AB$, it uniquely extends to a $K$-linear map $\psi\colon F\to K$. By non-degeneracy of the trace, there exists ${g}\in F$ such that $\psi(h)=\tr({g}h)$ for all $h\in F$.

    After replacing $Y$ by a suitable Zariski dense open affine subset $U$ and $X$ by its preimage, we can assume that $B$ is finitely generated and free as $A$-module, see for example \cite{stacks-project}*{Lemma 051R}. By further shrinking $U$, we can assume that $f$ is everywhere unramified and $g$ is regular on $X$. Then for $h\in B$ and $y\in Y(\R)$ we have
    \begin{equation}\label{eq:traceequality1}
        \sum_{x\in f^{-1}(y)\cap M}g_M(x)h(x)=G_h(y)=\tr({g}h)(y)=\sum_{x\in f^{-1}(y)}{g}(x)h(x).
    \end{equation}
    For the last equality see for example \cite{KummerDiss}*{Remark 2.2.8} or the proof of \cite{agostini}*{Lemma 4.16}. Since regular functions on affine varieties separate points, we can find for every $x\in M$ some $h_x\in B$ with $h_x(x)=1$ and which vanishes on every other preimage of $f(x)$. Applying (\ref{eq:traceequality1}) to this $h_x$ shows that $g_M(x)={g}(x)$. 
    
    For the second part, we choose $U$ to be an open dense affine subset of $Y$ such that $g$ does not vanish at any point of $f^{-1}(U)$ and agrees with $g_M$ on $M\cap f^{-1}(U)$. This is possible because $g$ is non-zero on a dense Zariski open subset of $M$ and $M$ is Zariski dense in $X$. 
    The equality in (\ref{eq:traceequality1}) reduces to
    \begin{equation}\label{eq:traceequality2}
        \sum_{x\in f^{-1}(y)\cap M}g(x)h(x)=\sum_{x\in f^{-1}(y)}{g}(x)h(x)
    \end{equation}
    for all $h\in B$ and $y\in Y(\R)$. Assume that there exists $x\in f^{-1}(U(\R))\smallsetminus M$. Then, as above, let $h_x\in B$ with $h_x(x)=1$ and which vanishes on every other preimage of $y=f(x)\in Y(\R)$. The left-hand side of (\ref{eq:traceequality2}) is zero because $x\not\in M$ and the right-hand side is $g(x)$. Thus $g(x)=0$ which contradicts our assumption that $g$ is non-zero on the preimage of $U$.
\end{proof}
\begin{lem}\label{lem:morphismfinite}
The morphism $f$ is finite. Equivalently, the ring $B$ is finitely generated as $A$-module.
\end{lem}
\begin{proof}
    Let $g$ be as in the previous lemma and choose a $K$-basis $x_1,\ldots,x_m\in B$ of $F$. Since $(x,y)\mapsto\tr(g xy)$ is nondegenerate, there exists a $K$-basis $y_1,\ldots,y_m$ of $F$ such that $\tr(gx_iy_j)=\delta_{ij}$ for all $i,j=1,\ldots,m$. Let $C\subset F$ be the $A$-module generated by $y_1,\ldots,y_m$. 
    Let $b\in B$ and write it as
    \begin{equation*}
     b=a_1y_1+\cdots+a_my_m
    \end{equation*}
    with $a_1,\ldots,a_m\in K$. Since $x_i\in B$ and $\tr(gB)\subset A$ by Lemma \ref{lem:tracealg}, it follows that 
    \begin{equation*}
    a_i=\tr(gx_ib)\in A    
    \end{equation*}
    for all $i=1,\ldots,m$ and thus $b\in C$. Hence $B\subset C$ which proves the claim because $A$ is noetherian and $C$ is finitely generated.
\end{proof}
\begin{proof}[Proof of Theorem \ref{thm:realfiberedcrit}]
    The morphism $f$ being finite was shown in Lemma \ref{lem:morphismfinite}. 
    Lemma \ref{lem:tracealg} shows that the degree of $f$ is the same as the degree of the covering $f|_M\colon M\to Y(\R)$. Because $f$ is finite and $Y$ smooth, no fiber of $f$ can have more elements than $\deg(f)$, with equality only if all points in the fiber are simple \cite{sga1}*{Expos\'e I, Th\'eor\`eme 10.11}. This shows that 
    \begin{equation*}
        X(\R)\subset f^{-1}(Y(\R))\subset M    \subset X(\R)
    \end{equation*}
    and that every point of $X(\R)$ is smooth. This proves the remaining claims.
\end{proof}
\subsection{Lighthouse manifolds are algebraic}
Let $M\subset\T^n$ be a closed submanifold of dimension $c$ that is a lighthouse with respect to $L\in\R^{n\times d}$ and $\nu$ the $L$-directed measure on $M$. By Proposition \ref{prop:lighthousesmooth} the manifold $M$ is smooth. We orient $M$ so that the density of $\nu$ is given by
\begin{equation*}
    \omega_L=\frac{1}{(2\pi i)^c}\sum_{I\in\binom{[n]}{c}}s(I)\cdot L_{ I}\cdot\frac{dz_{i_1}\wedge\cdots\wedge dz_{i_c}}{z_{i_1}\cdots z_{i_c}}
\end{equation*}
where in the sum the set $I$ consists of $i_1<\cdots<i_c$ and $s(I)=(-1)^{\sum_{i\in I}i-\sum_{i=1}^ci}$. The lighthouse condition implies that $\nu$ strongly avoids $\ker(L^t)$, i.e.~there exists $C>0$ such that
\begin{equation*}
    \|k\|\leq C\|L^tk\|
\end{equation*}
for all $k\in\supp(\hat\nu)$. 
Let $X\subset(\C^*)^n$ be the Zariski closure of $M$.
For a matrix $P\in\Z^{c\times n}$ we consider the map $$f=f_{P}\colon X\to(\C^*)^c,\, z\mapsto z^{P}.$$ The goal is to prove that $M$ is algebraic. More precisely, we will prove the following.
\begin{thm}\label{thm:lawtonalgebraic}
    Let $M\subset\T^n$ be a closed submanifold of dimension $c$ that is a lighthouse with respect to $L\in\R^{n\times d}$, $n=c+d$.
    If $X\subset(\C^*)^n$ is the Zariski closure of $M$, then $M=X\cap\T^n$ and every point of $M$ is a smooth point of $X$. Moreover, for every $P\in\Z^{c\times n}$ whose row-span is sufficiently close to the kernel of $L^t$, the map $f_P\colon X\to(\C^*)^c$ satisfies $f^{-1}_P(\T^c)=M$.
\end{thm}
If the row-span of $P$ is sufficiently close to the kernel of $L^t$ (in the sense of Lemma \ref{lem:lighthouselemma1}), then the restriction $M\to\T^c$ of $f$ is a covering map and the row-span of $P$ intersects the support of $\hat\nu$ only in the origin by Lemma \ref{lem:lighthouselemma1}. After multiplying a row of $P$ by $-1$ if necessary, this is an orientation preserving covering map. Thus there exists a smooth function $g_M$ on $M$ which is positive on $M$, and which satisfies $\omega_L=g_M\cdot f^*\eta$ where $\eta=\frac{1}{(2\pi i)^c}\frac{dw_1\wedge\cdots\wedge dw_c}{w_1\cdots w_c}$. Here the $w_i$ are the coordinate functions on $(\C^*)^c$. 
Finally, for $m\in\Z^n$ we define
\begin{equation*}
    G_m\colon \T^c\to\C,\, w\mapsto\sum_{z\in f^{-1}(w)\cap M}g_M(z)z^{-m}.
\end{equation*}
\begin{lem}
    For every $m\in\Z^n$ the function $G_m$ is a Laurent polynomial.
\end{lem}
\begin{proof}
    For $l\in\Z^c$ we calculate
    \begin{equation*}\label{eq:fourierFm}
        \hat{G}_m(l)=\int_{\T^c}G_m(w) w^{-l}\eta(w)=\int_Mg_M(z)z^{-P^tl-m}f^*\eta(z)=\hat\nu(P^tl+m).
    \end{equation*}
    By our choice of $P$ the measure $\nu$ strongly avoids the row-span of $P$, so Lemma \ref{lem:lighthouselemma2} shows that there are only finitely many $l\in\Z^c$ with $\hat\nu(P^tl+m)\neq0$. Thus $G_m$ is a Laurent polynomial.
\end{proof}
\begin{lem}
    The variety $X$ is of pure dimension $c$. The map $f=f_P\colon X\to(\C^*)^c$ is generically finite and dominant on every irreducible component of $X$. The degree of $f$ is the degree of the restriction $f|_M\colon M\to\T^c$.
\end{lem}
\begin{proof}
    Let $A$ and $B$ be the coordinate rings of $(\C^*)^c$ and $X$, respectively. Further, let $K$ be the rational function field of $(\C^*)^c$. Via $f$ we can consider $A$ as a subring of $B$ because $f$ is dominant. For $h\in B$ we consider the map
    \begin{equation*}
        G_h\colon\T^c\to\C,\, w\mapsto\sum_{z\in f^{-1}(w)\cap M}g_M(z)h(z)
    \end{equation*}
    which is, by the preceding lemma, a Laurent polynomial. The resulting map $\ell\colon B\to A,\,h\mapsto G_h$ is $A$-linear. We claim that for $0\neq b\in B$ there exists $h\in B$ with $\ell(hb)\neq0$. Indeed, for $0\neq b\in B$ there exists $x\in M$ with $b(x)\neq0$ and we can choose $h\in B$ which is nonzero at $x$ and zero at the other (finitely many) points of $f^{-1}(f(x))\cap M$. This implies in particular that $B$ is torsion-free as $A$-module. Indeed, let $0\neq b\in B$ and assume that $ab=0$ for $a\in A$. Then 
    \begin{equation*}
        0=\ell(abh)=a\cdot\underbrace{\ell(bh)}_{\neq0},
    \end{equation*}
    which shows that $a=0$ because $A$ does not have zero divisors.
    Hence $f$ is dominant on every irreducible component of $X$. Next, we extend $\ell$ to a $K$-linear map $$L:=B\otimes_AK\to K$$ which, by abuse of notation, we also denote by $\ell$. It remains to show that the dimension of $L$ as $K$-vector space is the degree $m$ of the covering $M\to\T^c$, so in particular finite, as this implies the remaining claims. Since the generic fiber cardinality of $f$ is equal to this dimension, it suffices to show that $\dim_K(L)\leq m$.
    Assume for the sake of a contradiction that there are $m+1$ $K$-linearly independent $b_0,\ldots,b_m\in B$. 
    Since the $K$-bilinear form $L\times L\to K$ defined by $(a,b)\mapsto\ell(ab)$ is nondegenerate, there are $h_0,\ldots,h_m\in B$ such that the matrix $(\ell(b_ih_j))_{0\leq i,j\leq m}$ has rank $m+1$. On the other hand, for every $w\in\T^c$, the matrix $(\ell(b_ih_j)(w))_{0\leq i,j\leq m}$ is, by definition of $\ell$, the sum of the $m$ rank one matrices $(g_M(z_k)b_i(z_k)h_j(z_k))_{0\leq i,j\leq m}$ where $z_1,\ldots,z_m$ are the elements of $f^{-1}(w)\cap M$. Thus the matrix $(\ell(b_ih_j))_{0\leq i,j\leq m}$ has rank at most $m$ over $K$, which is our desired contradiction.
\end{proof}
\begin{proof}[Proof of Theorem \ref{thm:lawtonalgebraic}]
    The preceding two lemmas show that the assumptions of Theorem \ref{thm:realfiberedcrit} are met with $Y=(\C^*)^c$. Here we choose the real structure on $(\C^*)^n$ and $(\C^*)^c$ given by the coordinate-wise involution $z\mapsto1/\bar{z}$.
\end{proof}

\section{The canonical algebraic models}\label{sec:final}
Recall that the \emph{amoeba} $\cA(X)\subset\R^n$ of a closed subvariety $X\subset(\C^*)^n$ is the image of $X$ under the map
\begin{equation*}
    \log|\cdot|\colon(\C^*)^n\to\R^n,\, (z_1,\ldots,z_n)\mapsto(\log|z_1|,\ldots,\log|z_n|).
\end{equation*}
\begin{Def}\label{def:verynicedef}
    Let $X\subset(\C^*)^n$ be a closed subvariety of pure dimension $c=n-d$ and $V\in\Gr(d,n)$. Consider the following conditions:
    \begin{enumerate}
    \item The variety $X$ is invariant under the coordinate-wise involution $z\mapsto1/\bar{z}$.
     \item There exists an open neighborhood $U\subset\Gr(d,n)$ of $V$ such that
    \begin{equation*}
        V'\cap\cA(X)\subset\{0\}
    \end{equation*}
    for all $V'\in U$.
    \item Every point of $X\cap\T^n$ is a smooth point of $X$.
    \end{enumerate}
 We say that $X$ is \emph{essentially Lee--Yang} with respect to $V$ if conditions (1) and (2) are satisfied. If, in addition, condition (3) holds, we say that $X$ is \emph{\varietyname{}} with respect to $V$.
 We say that $X$ is \emph{essentially (strict) Lee--Yang} with respect to a matrix $L\in\R^{n\times d}$ of rank $d$ if it has this property with respect to the range of $L$.
\end{Def}
We will prove that \varietyname{} varieties are exactly the algebro-geometric counterparts of Delone Fourier quasicrystals.
\begin{lem}\label{lem:verynicedelone}
    Let $L\in\R^{n\times d}$ have rank $d$, let $X\subset(\C^*)^n$ be \varietyname{} with respect to $L$ and let $M=X\cap\T^n$. Furthermore, let $Y=\{z\in\C^n\mid\exp(2\pi i z)\in X\}$ and $Y(\R)=Y\cap\R^n$. Then the following holds:
    \begin{enumerate}
        \item $X\cap\T^n$ is Zariski dense in $X$.
        \item The tangent space of $Y(\R)$ at every point is transverse to the range of $L$.
        \item The set $\Lambda(X,L)$ is Delone.
    \end{enumerate}
\end{lem}
\begin{proof}
    Let $U\subset\Gr(d,n)$ be the open set from part (2) of Definition \ref{def:verynicedef}.
    Let $P\in\Z^{c\times n}$ be a matrix whose kernel is in $U$ such that under the map $$f_{P}\colon X\to(\C^*)^c,\, z\mapsto z^{P}$$ every irreducible component of $X$ dominates $(\C^*)^c$ (which exists for example by \cite{tropical}*{Corollary 3.2.13, Theorem 3.3.8}). 
    This implies that there is a dense subset of $\T^c$ in the image of each irreducible component $X_0$ of $X$ under $f_P$.
    On the other hand, if $x\in X_0$ with $f_P(x)\in\T^c$, then $\log|x|\in\ker(P)$. Now by condition (2) of Definition \ref{def:verynicedef}, we must have $\log|x|=0$, meaning that $x\in X_0\cap\T^n$.
    Thus $X_0\cap\T^n$ is mapped by $f_P$ to a dense subset of $\T^c$. Thus it must have real dimension at least $c$, which shows that it is Zariski dense in $X_0$.
    This shows part (1) of the claim.

    For part (2) we follow the same strategy as in the proof of \cite{alon2025higher}*{Lemma~4.17}. 
    We first note that, because $X$ is invariant under the coordinate-wise involution $z\mapsto1/\bar{z}$, it follows that $Y$ is invariant under complex conjugation. This implies that the tangent space $T_xY$ of $Y$ at every real point $x\in Y(\R)$ is the $\C$-span of the tangent space of $Y(\R)$.
    We also note that $\cA(X)$ is equal to the set
    \begin{equation*}
        \cA(X)=\{-2\pi\cdot\Im(z)\mid z\in Y\}.
    \end{equation*}
    This shows that for $z\in Y$ and $V\in U$ we can have $\Im(z)\in V$ only if $z\in Y(\R)$.
    Now for every $0\neq w\in T_xY(\R)$ we have to show that $w$ is not in the range of $L$. Then $\xi=i\cdot w\in T_xY$ and we choose a smooth path $\gamma\colon[0,1]\to Y$ with $\gamma(0)=x$ and $\gamma'(0)=\xi$. Then for small $t\in(0,1]$ we have $\Im(\gamma(t))\neq0$ and thus 
    \begin{equation*}
        \Im(\frac{1}{t}(\gamma(t)-\gamma(0)))=\frac{1}{t}\Im(\gamma(t))\not\in V
    \end{equation*}
    for any $V\in U$. Since $U$ is an open neighborhood of the range of $L$, this shows that the limit $w=\Im(\gamma'(0))$ is not in the range of $L$.

    For proving part (3) we consider the submanifold $\Sigma=Y(\R)/\Z^n$ of $(\R/\Z)^n$. Then
    \begin{equation*}
        \Lambda:=\Lambda(X,L)=\{x\in\R^d\mid \overline{Lx}\in\Sigma\}.
    \end{equation*}
    We first prove that $\Lambda$ is uniformly discrete. Assume for the sake of a contradiction that there is a sequence of distinct $x_j,y_j\in\Lambda$ such that $s_j:=x_j-y_j$ tends to zero. Since $\Sigma$ is compact, we can assume, after passing to subsequences, that the sequences $(\overline{Lx_j})_j$ and $(\overline{Ly_j})_j$ both converge to the same point $x\in\Sigma$ and that $\frac{s_j}{\|s_j\|}$ converges to some $0\neq s\in\R^d$. But then $Ls$ is in the tangent space of $\Sigma$ at $x$ as the limit of secants of pairs of points converging to $x$. This contradicts part (2). Finally, we prove that $\Lambda$ is relatively dense. Consider the map
    \begin{equation*}
        F\colon\Sigma\times\R^d\to(\R/\Z)^n,\, (a,b)\mapsto a+\overline{Lb}
    \end{equation*}
    which is locally a diffeomorphism by part (2), so in particular open. We claim that it is surjective. Indeed, if $H$ is the closure of the image of $x\mapsto\overline{Lx}$, then the image of $F$ is $\Sigma+H$. Thus $\Sigma+H$ is closed and open, which implies surjectivity of $F$ because $(\R/\Z)^n$ is connected.
    In particular, the family $F(\Sigma\times B_R)$, $R>0$, of open subsets covers $(\R/\Z)^n$. By compactness we find $R>0$ with $F(\Sigma\times B_R)=(\R/\Z)^n$. This shows that, for every $y\in\R^d$ there exists $v\in\R^d$ with $\|v\|<R$ and $x\in\Sigma$ such that $\overline{Ly}=\overline{Lv}+x$, so $y-v\in B_R(y)\cap\Lambda$.
\end{proof}

\begin{lem}\label{lem:veryniceislighthouse}
 Let $L\in\R^{n\times d}$ be a matrix of rank $d$. Let $X\subset(\C^*)^n$ be \varietyname{} with respect to $L$. Then $M=X\cap\T^n$ is a lighthouse manifold with respect to $L$.
\end{lem}

\begin{proof}
    By Lemma \ref{lem:verynicedelone} the manifold $M$ is nonempty and transverse to $L$. Let $U$ be the neighborhood of the range of $L$ from part (2) of Definition \ref{def:verynicedef}. By \cite{alon2025higher}*{Lemma 8.1} there exists a non-singular $n\times n$ matrix $S$ with integer entries such that $SL$ has only positive $d\times d$ minors and such that the range of $S^{-1}L'$ is in $U$ for every matrix $L'\in\R^{n\times d}$ of rank $d$ with only nonnegative $d\times d$ minors. Let $m\in\N$ such that $A=mS^{-1}$ has integer entries and consider the morphism
    \begin{equation*}
        f_A\colon(\C^*)^n\to(\C^*)^n,\, z\mapsto z^A.
    \end{equation*}
    By decomposing $A$ into Smith normal form, one can see that $f_A$ is finite and unramified and satisfies $f_A^{-1}(\T^n)=\T^n$. Thus, letting $X'=f_A^{-1}(X)$, we have that $M'=X'\cap\T^n$ is contained in the smooth part of $X'$. Further note that all $d\times d$ minors of  $\tilde{L}=A^{-1}L$ are positive. We have $\cA(X')=A^{-1}\cA(X)$, so that every linear subspace in $\Gr_{\geq}(d,n)\subset A^{-1}U$ intersects $\cA(X')$ only in the origin. Thus by \cite{alon2025higher}*{Lemma 8.3} it follows that $X'$  is a strict Lee--Yang variety. Therefore, by Example \ref{ex:leeyanglighthouse} the manifold $M'$ is a lighthouse with respect to $\tilde{L}$. Now if $\nu'$ is the $\tilde{L}$-directed measure on $M'$, then $\nu:=(f_A)_*\nu'$ is the $L$-directed measure on $M$. Hence $\hat{\nu}(k)=\hat\nu'(A^tk)$ and the lighthouse property of $M'$ implies the lighthouse property of $M$.
\end{proof}
\begin{lem}\label{lem:coordchange}
    Let $L\in\R^{n\times d}$ have rank $d$ and let $X\subset(\C^*)^n$ be \varietyname{} with respect to $L$.
    Then there is a strict Lee--Yang variety $\tilde{X}\subset(\C^*)^n$ and an $n\times d$ matrix $\tilde{L}$ whose range is in $\Gr_{+}(d,n)$ such that 
    $\Lambda(X,L)=\Lambda(\tilde{X},\tilde{L})$. In particular, the set $\Lambda(X,L)$ is a Delone Fourier quasicrystal.
\end{lem}
\begin{proof}
    The claim follows from \cite{alon2025higher}*{Theorem 8.4} because $\Lambda(X,L)$ is Delone by Lemma \ref{lem:verynicedelone}.
\end{proof}

\begin{lem}\label{lem:rationalimpliesamoeba}
    Let $X\subset(\C^*)^n$ be a closed subvariety of pure dimension $c=n-d$. Assume that there is an open subset $U\subset{\Gr}(d,n)$ such that $\cA(X)\cap V=\{0\}$ for every rational $V\in U$. Then $\cA(X)\cap V=\{0\}$ for every $V\in U$.
\end{lem}
\begin{proof}
    Without loss of generality, we can assume that $X$ is irreducible. Assume for the sake of a contradiction that there exists $0\neq a\in\cA(X)$ which is contained in some $V\in U$. Every $a'$ in some open neighborhood $W$ of $a$ is contained in some $V'\in U$. The preimage of $W$ under $\log|\cdot|$ is a nonempty open subset of $X$ so that it contains a smooth point of $X$. Therefore, without loss of generality, we can assume that for the above $a$ there exists a smooth point $z\in X$ with $a=\log|z|$. We first consider the case when the differential of $\log|\cdot|\colon X\to\R^n$ at $z$ has rank strictly larger than $c$. After replacing $V$ by a nearby element of $U$ if necessary, we can assume that $V$ and the image of this differential at $z$ intersect transversely. Then every $V'$ in a neighborhood of $V$ also intersects $\cA(X)$ in a nonzero point. Since such a neighborhood also contains rational planes, we arrive at a contradiction. Thus the rank of the differential of $\log|\cdot|$ at $z$ must be at most $c$. Since we chose $z$ to be an arbitrary smooth point of $X$ in the preimage of $W$, the differential at every smooth point of $X$ in this open preimage, and thus on every smooth point of $X$, is at most $c$. This shows that $\cA(X)$ has dimension at most $c$. By \cite{nissesottile}*{Theorem 4.6} this implies that $X$ is the coset of an algebraic subtorus of $(\C^*)^n$ of dimension $c$. Thus the amoeba $\cA(X)$ is a rational affine subspace of $\R^n$ of dimension $c$. Let $V'\in U$ be some rational subspace that is complementary to the vector space of pairwise differences of elements of $\cA(X)$, so that it intersects $\cA(X)$. By assumption this intersection point must be the origin. Thus $\cA(X)$ is a rational linear subspace of $\R^n$. Now let $0\neq a_1\in\cA(X)\cap V$ and $a_1,\ldots,a_d$ a basis of $V$. If we replace $a_1$ by a close enough rational non-zero element of $\cA(X)$ and $a_2,\ldots,a_d$ by close enough rational vectors, then their span is a rational space in $U$ that intersects $\cA(X)$ non-trivially. This is a contradiction to our assumption.
\end{proof}

\begin{cor}\label{cor:lighthouseisverynice}
    Let $M\subset\T^n$ be a closed submanifold of dimension $c$ that is a lighthouse with respect to $L\in\R^{n\times d}$, $n=c+d$.
    The Zariski closure $X\subset(\C^*)^n$ of $M$ is \varietyname{} with respect to $L$ and $M=X\cap\T^n$.
\end{cor}

\begin{proof}
    By Theorem \ref{thm:lawtonalgebraic} we have $M=X\cap\T^n$ and condition (3) of Definition \ref{def:verynicedef}. Condition (1) of Definition \ref{def:verynicedef} is automatically satisfied since $X$ is the Zariski closure of a subset of $\T^n$.

    For proving part (2) of Definition \ref{def:verynicedef} let $U\subset {\Gr}(d,n)$ be an open neighborhood of the column span $V$ of $L$ such that the orthogonal complement of every $V'\in U$ is sufficiently close to the kernel of $L^t$ in the sense of Theorem \ref{thm:lawtonalgebraic}. Thus if $V'\in U$ is rational, then it is the kernel of some matrix $P\in\Z^{c\times n}$ for which the map $f_P\colon X\to(\C^*)^c$ satisfies $f^{-1}_P(\T^c)=X\cap\T^n$. Now let $a\in\cA(X)\cap V'$ and let $z\in X$ such that $\log|z|=a$. Because $a$ is in the kernel of $P$, we have
    \begin{equation*}
        \log|f_P(z)|=P\log|z|=Pa=0,
    \end{equation*}
    which implies that $f_P(z)\in\T^c$. Because $f^{-1}_P(\T^c)=X\cap\T^n$, we get $z\in\T^n$ and hence $a=0$. Thus Lemma \ref{lem:rationalimpliesamoeba} implies condition (2) of Definition \ref{def:verynicedef}. 
\end{proof}

\begin{thm}
    Let $\emptyset\neq M\subset\T^n$ be a closed submanifold of dimension $c=n-d$. Assume that $M$ is transverse with respect to the matrix $L\in\R^{n\times d}$ of rank $d$ whose rows are linearly independent over $\Q$. Then the following are equivalent:
    \begin{enumerate}
        \item $M$ is lighthouse with respect to $L$.
        \item $M=\T^n\cap X$ for a closed subvariety $X\subset(\C^*)^n$ that is \varietyname{} with respect to $L$.
        \item $\Lambda=\Lambda(M,L)$ satisfies Assumption~\ref{ass:crystalline} {with spectrum contained in the group generated by the rows of $L$}.
    \end{enumerate}%
\end{thm}

\begin{proof}
    For the direction $(1)\Rightarrow(2)$ we can take as $X$ the Zariski closure of $M$ in $(\C^*)^n$ by Corollary \ref{cor:lighthouseisverynice}. The direction $(2)\Rightarrow(1)$ is Lemma \ref{lem:veryniceislighthouse} and the direction $(1)\Rightarrow(3)$ is Proposition \ref{prop:lighthousegivesFQ}. Finally, we prove the direction $(3)\Rightarrow(1)$. Let $\Gamma$ be the group generated by the spectrum of $\Lambda$, and let $m=\rank(\Gamma)\le n$. By assumption, there exists $A\in\Z^{m\times n}$ of full rank such that the rows of $L_0:=AL$ form a $\Z$-basis of $\Gamma$. Consider the map $$f_A\colon\T^n\to\T^m,\, z\mapsto z^A.$$ By transversality of $M$ the image $M_0$ of $M$ under $f_A$ is Lawton's manifold of $\Lambda$. If $\nu_0$ is the $L_0$-directed measure on $M_0$ and $\nu$ the $L$-directed measure on $M$, then $\supp(\hat\nu)=A^t\supp(\hat\nu_0)$. Thus, since $M_0$ is lighthouse with respect to $L_0$ by Corollary \ref{cor:lawtonlighthouseandsmooth}, it follows that $M$ is a lighthouse with respect to $L$.
\end{proof}

\begin{lem}\label{lem:lawtonverynice}
    Let $\Lambda\subset\R^d$ satisfy Assumption~\ref{ass:crystalline}. Let $L\in\R^{n\times d}$ be a matrix whose rows form a $\Z$-basis of the group generated by the spectrum of $\Lambda$. Finally, let $M\subset\T^n$ be Lawton's manifold, as defined in Section \ref{sec:lawton}, and $X\subset(\C^*)^n$ its Zariski closure. Then the following holds:
    \begin{enumerate}
     \item $\Lambda=\Lambda(M,L)=\Lambda(X,L)$.
     \item $X$ is \varietyname{} with respect to $L$.
     \item If $z_0\in\T^n$ such that $z_0\cdot M=M$, then $z_0=(1,\ldots,1)$.
    \end{enumerate}
\end{lem}
\begin{proof}
 Part (1) of the claim is covered by Theorem \ref{thm: Lawton} and Corollary \ref{cor:lighthouseisverynice}.

 Part (2) follows from Corollary \ref{cor:lighthouseisverynice} and Corollary \ref{cor:lawtonlighthouseandsmooth}.

    Finally, for part (3), let $z_0\in\T^n$ such that $z_0\cdot M=M$. Let $\nu=\nu_L$ be the $L$-directed measure on $M$. We note that the pushforward of $\nu$ under multiplication by $z_0$ equals $\nu$. Therefore, for every $k\in\Z^n$ with $z_0^k\neq1$, we have
    \begin{equation*}
     \hat\nu(k)=\int z^{-k}d\nu(z)=\int z_0^{-k}z^{-k}d\nu(z)=z_0^{-k}\hat\nu(k)\Rightarrow\hat\nu(k)=0.
    \end{equation*}
    Thus $\supp(\hat\nu)$ is contained in the subgroup $G$ of $\Z^n$ of all $k$ with $z_0^k=1$. Theorem \ref{thm: Lawton} then implies that the spectrum $S(\Lambda)$ is contained in $L^tG$. However, by assumption, the group generated by $S(\Lambda)$ is $L^t\Z^n$. By injectivity of $L^t$, this implies that $G=\Z^n$ and thus $z_0^k=1$ for all $k\in\Z^n$. This shows that $z_0=(1,\ldots,1)$.
\end{proof}

\begin{proof}[Proof of Theorem \ref{thm:main}]
    The claim follows directly from Lemma \ref{lem:lawtonverynice} and Lemma \ref{lem:coordchange}.
\end{proof}
We say that $X\subset(\C^*)^n$ has no torus stabilizer if the only $z_0\in\T^n$ such that $z_0\cdot X=X$ is the trivial element $z_0=(1,\ldots,1)$.
\begin{thm}\label{thm:bijection1}
    Let $L\in\R^{n\times d}$ be a matrix of rank $d$ whose rows are linearly independent over $\Q$ and $\Gamma\subset\R^d$ be the group generated by the rows of $L$. 
    We consider the following three sets:
    \begin{align*}
    \textnormal{FQ}(L)&:=\{\Lambda\subset\R^d\mid \Lambda\textrm{ satisfies Assumption~\ref{ass:crystalline} and its spectrum generates }\Gamma\},\\
    \textnormal{LH}(L)&:=\{\textrm{Lighthouse manifolds }M\subset\T^n\textrm{ w.r.t. }L\textrm{ without torus stabilizer}\},\\
    \textnormal{LY}(L)&:=\{\textrm{ess. strict Lee--Yang varieties }X\subset(\C^*)^n\textrm{ w.r.t.}\,L\textrm{ w/o torus stabilizer}\}.
    \end{align*}
    The maps $\textnormal{LY}(L)\to\textnormal{LH}(L),\, X\mapsto X\cap\T^n$ and $\textnormal{LH}(L)\to\textnormal{FQ}(L),\, M\mapsto\Lambda(M,L)$ are bijections. Their inverse maps are given by 
    \begin{equation*}
        M\mapsto \overline{M}^{\rm Zar} \textrm{ and } \Lambda\mapsto\overline{\{\exp(2\pi i Lx)\mid x\in\Lambda\}}.
    \end{equation*}
\end{thm}
\begin{proof}
    We first show that the map $\textnormal{LY}(L)\to\textnormal{LH}(L),\, X\mapsto X\cap\T^n$ is well-defined.
    By Lemma \ref{lem:veryniceislighthouse} we only have to show that there is no torus stabilizer. But if $X\cap\T^n$ had a torus stabilizer, then $X$ would have one too by part (1) of Lemma \ref{lem:verynicedelone}. The map $M\mapsto \overline{M}^{\rm Zar}$ is indeed the inverse by Corollary \ref{cor:lighthouseisverynice} and part (1) of Lemma \ref{lem:verynicedelone}.
    
    We next show that the map $\textnormal{LH}(L)\to\textnormal{FQ}(L),\, M\mapsto\Lambda(M,L)$ is well-defined. By Proposition \ref{prop:lighthousegivesFQ} $\Lambda(M,L)$ is a Delone Fourier quasicrystal and thus satisfies Assumption \ref{ass:crystalline}. It remains to show that the spectrum of $\Lambda(M,L)$ generates $\Gamma$. But this follows from the description of the spectrum of $\Lambda(M,L)$ in Proposition \ref{prop:lighthousegivesFQ} and the fact that the Fourier support of the $L$-directed measure $\nu_L$ on $M$ generates $\Z^n$ (otherwise $\nu_L$ and thus $M$ would have a torus stabilizer).  The map that assigns Lawton's manifold to $\Lambda\in\textnormal{FQ}(L)$ is indeed the inverse by Corollary \ref{cor:lawtonlighthouseandsmooth}, Lemma \ref{lem:lawtonverynice}, and Proposition \ref{prop:lighthousegivesFQ}.
\end{proof}
\begin{cor}
    Let $M\subset\T^n$ be a closed submanifold of dimension $c=n-d$ and $L\in\R^{n\times d}$ a matrix of rank $d$ whose rows are linearly independent over $\Q$. Then the following are equivalent:
    \begin{enumerate}
        \item $M$ is lighthouse with respect to $L$ and has no torus stabilizer.
        \item $M=\T^n\cap X$ for a closed subvariety $X\subset(\C^*)^n$ that is \varietyname{} with respect to $L$ and which has no torus stabilizer.
        \item $\Lambda=\Lambda(M,L)$ satisfies Assumption~\ref{ass:crystalline} and $M$ is Lawton's manifold of $\Lambda$.
    \end{enumerate}
\end{cor}
\begin{proof}
    This follows directly from the previous theorem.
\end{proof}

\section*{Acknowledgements}

We thank {Henry Cohn, Larry Guth,
David Jerison, Pavel Kurasov, Wayne Lawton, Peter Sarnak, and Cynthia Vinzant}
for useful conversations.
{We further thank Kaie Kubjas and Bernd Sturmfels for organizing the learning seminar \emph{Holonomic D-cember} at UC Berkeley in 2014 where the second author learned about holonomic $D$-modules and the $b$-function.}
The first author was supported by the Simons Foundation
through the Collaboration on Localization of Waves
(Grant No.~1027116, D.J.).

\section*{Use of AI and concurrent work}

We consulted {different AI systems} at various stages of this work, {mainly
for assistance with wording and exposition, proofreading, and literature research.} The authors developed the proofs
and {wrote} the manuscript themselves, over a period of two years, and take full responsibility for its content.

We are aware of the {completely AI generated} preprint \emph{Fourier quasicrystals with
positive integer masses are sections of Lee--Yang cycles}
by Dongsheng Wei, { first submitted to the arXiv} on October 7, 2026\footnote{ arXiv Identifier	2610.10652}.
The present work was developed independently and in parallel.
We have not verified the proofs in that preprint and make
no claim regarding their correctness. We have not assessed
the extent of overlap between the two manuscripts.

\appendix
\section{A direct proof of Lawton's manifold theorem}\label{app:lawton}
For $N\in\N$ and $h\in\mathcal S(\R^d)$, define the
weighted supremum norm
\[
    \|h\|_N
    :=
    \sup_{x\in\R^d}(1+\|x\|)^N|h(x)|.
\]
It is useful for estimates of Schwartz functions over annuli. For any $h\in\mathcal{S}(\R^d)$, any $N\in\N$ and any $R>1$,
\begin{equation}\label{eq: Am estimates}
    |h(x)|\le\frac{\|h\|_N}{(1+2^mR)^{N}},\quad \forall x\in A_m(R),
        \qquad m\in\N_0, 
\end{equation}
\[
        A_m(R):=
        \{x\in\R^d:2^mR\le\|x\|<2^{m+1}R\}.
    \]
\begin{lem}\label{lem:vagueSchwartz}
    Let $\delta,A>0$, and consider a sequence of measures
    \[
        \mu_j=\sum_{x\in\Lambda_j}a_j(x)\delta_x,
    \]
    where each $\Lambda_j\subset\R^d$ is
    $\delta$-uniformly discrete and $|a_j(x)|\le A$
    for all $j\in\N$ and $x\in\Lambda_j$.
    These measures are tempered, and the following
    conditions are equivalent:
    \begin{enumerate}
        \item The sequence $(\mu_j)_{j=1}^\infty$ converges vaguely.
        \item The sequence $(\mu_j)_{j=1}^\infty$ converges in
        $\mathcal S'(\R^d)$.
        \item The sequence $(\widehat{\mu_j})_{j=1}^\infty$ converges
        in $\mathcal S'(\R^d)$.
    \end{enumerate}
    In this case, denoting the vague limit by $\mu$,
    the distributional limits in (2) and (3) are
    $\mu$ and $\widehat\mu$, respectively.
    Moreover, whenever these limits exist, for every
    integer $N>d$ there is a constant
    $C_{\delta,A,d,N}>0$ such that
    \[
        \sup_j\int_{\R^d}|h|\,d|\mu_j|
        +\int_{\R^d}|h|\,d|\mu|
        \le C_{\delta,A,d,N}\|h\|_N
        \qquad
        \text{for all $h\in\mathcal S(\R^d)$.}
    \]
     If, in addition, each $\widehat{\mu_j}$ is a
    pure-point measure
    \[\widehat{\mu_j}=\sum_{\xi\in S}c_{j,\xi}\delta_\xi,\]
    supported in a common locally
    finite set $S\subset\R^d$, then these conditions
    are equivalent to the point-wise convergence $c_{j,\xi}\to b_{\xi}$ to some $b_\xi$ for every $\xi\in S$.
    In that case,
    \[
        \widehat\mu
        =
        \sum_{\xi\in S}b_\xi\delta_\xi.
    \]
\end{lem}

\begin{proof}
    For any $R>0$, the $\delta/2$-balls centered around the points of $\Lambda_j\cap B_{R}$ are disjoint and contained in $B_{R'}$ with $R'=R+\frac{\delta}{2}$, hence the number of points in $\Lambda_j\cap B_{R}$ is bounded by $\mathrm{vol}(B_{R'})/\mathrm{vol}(B_{\frac{\delta}{2}})=\left(1+\frac{2R}{\delta}\right)^d$, so
    \[
        |\mu_j|(B_R)
        \le A\,\#(\Lambda_j\cap B_{R}(0))
        \le A\left(1+\frac{2R}{\delta}\right)^d
        \le C_0(1+R)^d,
    \]
    where $C_0$ depends only on $\delta,A,d$.

    These uniform bounds on total variation imply,
    by weak-* compactness and a diagonal argument,
    that every subsequence has a vaguely convergent
    further subsequence. By taking the supremum of $\int gd|\mu_j|$ over all $j$ and all
    $g\in C_c(B_R)$ with $0\le g\le 1$,  we conclude that any vague sub-sequential limit $\mu$
    satisfies the same growth bound. This proves temperedness of all measures in the sequence and all possible sub-sequential limits. 

    We next estimate the tails uniformly. Fix an
    integer $N>d$ and $R\ge1$, and let $\eta$ be any measure among the $\mu_j$ and
    their vague limits. Then, \eqref{eq: Am estimates} gives
    \begin{align*}
        \int_{\|x\|\ge R}|h(x)|\,d|\eta|(x)
        &=
        \sum_{m=0}^\infty
        \int_{A_m(R)}|h(x)|\,d|\eta|(x)\\
        &\le
        \|h\|_N
        \sum_{m=0}^\infty
        (1+2^mR)^{-N}
        |\eta|(B_{2^{m+1}R}(0))\\
        &\le
        C_0\|h\|_N
        \sum_{m=0}^\infty
        \frac{(1+2^{m+1}R)^d}{(1+2^mR)^N}\\
        &\le
        C_0\,2^{2d}\|h\|_N R^{d-N}
        \sum_{m=0}^\infty 2^{-m(N-d)}\\
        &=
        K_N\|h\|_N R^{d-N},
    \end{align*}
    where we used
    $1+2^{m+1}R\le2^{m+2}R$, 
    $1+2^mR\ge2^mR$, and $K_N:=\frac{C_0\,2^{2d}}{1-2^{d-N}}$. In particular, since $\sup_{\|x\|<1}|h(x)|\le\|h\|_{N}$, then substituting $R=1$ gives 
    \[\int_{\R^d}|h(x)|\,d|\eta|(x)\le C_0 2^d\|h\|_N+ K_N\|h\|_N, \]
    and setting $C_{\delta,A,d,N}=2(C_{0} 2^d+K_{N})$ gives the needed upper bound. 

    Suppose now that $\mu_j\to\mu$ vaguely.
    Choose a smooth cut-off $\chi_R\in C_c^\infty(\R^d)$ with
    $0\le\chi_R\le1$ and $\chi_R=1$ on $B_R$.
    Then
    \[
        |\langle\mu_j-\mu,h\rangle|
        \le
        |\langle\mu_j-\mu,\chi_Rh\rangle|
        +2K_N\|h\|_N R^{d-N}.
    \]
    For fixed $R$, the first term tends to zero
    by vague convergence. Letting $R\to\infty$
    proves $\mu_j\to\mu$ in $\mathcal S'$.

    Conversely, if $(\mu_j)$ converges in
    $\mathcal S'$, every vaguely convergent
    subsequence has the same distributional limit
    by the preceding argument. 
    This proves the equivalence of (1) and (2).
    The equivalence of (2) and (3) follows because
    the Fourier transform is an automorphism of
    $\mathcal S'$.

    Finally, suppose that the $\widehat{\mu_j}$ are
    pure-point measures supported in a common locally
    finite set $S$. Every $\xi\in S$ has a neighborhood $U_\xi$ disjoint from $S\setminus\{\xi\}$, and we can choose a function $f_{\xi}\in C_{c}^{\infty}(\R^d)$ supported in $U_{\xi}$ with $f_{\xi}(\xi)=1$. This means $\langle\widehat{\mu_{j}},f_{\xi}\rangle=c_{j,\xi}$ for all $j$ and all $\xi\in S$, and vague or $\mathcal{S}'$ convergence implies the convergence of these Fourier coefficients.

    Conversely, suppose that $\lim_{j\to\infty}c_{j,\xi}=b_{\xi}$ for every $\xi\in S$. As discussed previously, every subsequence of $(\mu_j)_{j=1}^{\infty}$ has a converging subsequence, so our goal is to show that all these sub-sequential limits agree. Let $\mu_*:=\lim_{k\to\infty}{\mu_{j_k}}$
    be any such sub-sequential limit, then it is both a vague limit and an $\mathcal{S}'$ limit. For any $g\in C_{c}^{\infty}(\R^d)$
    \[
        \langle\mu_*,\widehat{g}\rangle=
        \lim_{k\to\infty}\sum_{\xi\in S}c_{j_k,\xi}
        g(\xi)
        =
        \sum_{\xi\in S}b_\xi g(\xi),
    \]
    The sums are finite because $S$ is locally finite, so we can exchange summation and limit.
 Hence, every sub-sequential limit $\mu_*$ has the same Fourier transform, and are therefore equal. The relative compactness allows to conclude that there is a unique limit $\mu$ with the needed Fourier transform.
\end{proof}
Here is a proof of Theorem~\ref{thm: Lawton}:
\begin{proof}[Proof of Theorem~\ref{thm: Lawton}]
    Write $S=S(\Lambda)$ and
    $\varphi(y)=\exp(2\pi iLy)$, and put $c_\xi=0$
    for $\xi\notin S$. Let $\Gamma=\spann_{\Z}(S(\Lambda))$ and $n=\rank(\Gamma)$, which is finite by our assumption. Let $L\in\R^{n\times d}$ whose rows form
    a $\Z$-basis of $\Gamma$. By definition, $L^t:\Z^n\to\Gamma$ is an isomorphism, and as a result, $\varphi(\R^d)$ is dense
    in $\T^n$. To see why, notice that $\overline{\varphi(\R^d)}=\T^n$ if and only if $\overline{L\R^d/\Z^n}=\R^n/\Z^n$. The orthogonal complement to the tangent space of $\overline{L\R^d/\Z^n}$ is the span of $\ker(L^t)\cap\Z^n$. Thus $\ker(L^t)\cap\Z^n=\{0\}$ implies that $\overline{\varphi(\R^d)}=\T^n$.
    We also claim that $L$ has full rank. Otherwise,
    a nonzero vector $v\in\ker(L)$ would make
    $\widehat{\mu_\Lambda}$ invariant to $\widehat{\mu_\Lambda}(\{\xi\})\mapsto \widehat{\mu_\Lambda}(\{\xi\})e^{2\pi i t \langle \xi,v\rangle}$ for any $t\in\R$. This would imply
    $\Lambda-tv=\Lambda$ for every $t$, contradicting the assumption that $\Lambda$ is non-empty and uniformly discrete.

Define $\iota:\Gamma\to\Z^n$ to be the inverse of $L^t$. Denote  
\[\Omega=\{k\in\Z^n\ :\ L^tk\in S\},\quad \text{and$\quad b_{k}=c_{L^tk}\quad$ for every $k\in\Z^n$}.\]

    For $z\in\T^n$, choose $y_j\in\R^d$ with
    $\varphi(y_j)\to z$. The translated counting
    measures satisfy
    \[
        \widehat{\mu_{\Lambda-y_j}}
        =\sum_{\xi\in S}c_{\xi}e^{2\pi i \langle y_j,\xi\rangle }\delta_{\xi}=
        \sum_{k\in\Omega}
        b_k \varphi(y_j)^k\delta_{L^tk},
    \]
    and since $\varphi(y_j)^k\to z^k$ for all $k$, there is a unique limit $\mu_z$, with
    \begin{equation}\label{eq:torusfamilyFourier}
        \widehat{\mu_z}
        =
        \sum_{k\in \Omega}
        b_{k}z^k\delta_{L^tk}.
    \end{equation}
    This determines $\mu_z$, independently
    of the chosen sequence. A vague limit of
    counting measures of $\delta$-uniformly discrete sets must be a counting measure of a $\delta$-uniformly discrete set. Define $\Lambda_z$ to be the set such that
    $\mu_z=\mu_{\Lambda_z}$.

Comparing Fourier coefficients gives
\begin{equation}\label{eq:torusfamilytranslation}
        \Lambda_{\mathbf{1}}=\Lambda,
        \qquad
        \Lambda_{z\varphi(y)}=\Lambda_z-y.
    \end{equation}

    For any converging sequence in the torus $z_j\to z$, the coefficients converge $b_kz_j^k\to b_k z^k$, for all $k\in\Z^n$. Applying the coefficient-convergence criterion in Lemma~\ref{lem:vagueSchwartz} gives a vague and $\mathcal{S}'$ limit $\mu_{z_j}\to\mu_z$. In particular, $z\mapsto \int fd\mu_z$ is continuous for any fixed $f\in\mathcal{S}(\R^d)$. Given $f\in\mathcal{S}(\R^d)$, we define the continuous
    functions
    \[
        \psi_f(z):=\int \widehat{f}d\mu_z
        =
        \sum_{x\in\Lambda_z}\widehat f(x).
    \]
    Equation \eqref{eq:torusfamilyFourier} suggests that the $k$-th Fourier coefficient of $\psi_f$ should be $b_{k}f(L^tk)$. To justify this, we need a cut-off convergence argument. 
    Choose $\chi\in C_c^\infty(\R^d)$ equal to one
    near zero and put $f_R(\xi)=\chi(\xi/R)f(\xi)$.
    By \eqref{eq:torusfamilyFourier},
    \[
        \psi_{f_R}(z)
        =
        \sum_{k\in \Omega}
        b_k z^k f_R(L^tk),
    \]
    where the sum is only over $k\in\Omega$ with  $L^tk\in S\cap\supp(f_R) $, so it is finite. Since $f_R\to f$
    in $\mathcal S$, the uniform bound in
    Lemma~\ref{lem:vagueSchwartz} gives, for $N>d$,
    \[
        \sup_{z\in\T^n}|\psi_{f_R}(z)-\psi_f(z)|
        \le C_N\|\widehat{f_R-f}\|_N
        \rightarrow 0,\quad \text{as $R\to\infty$}.
    \]
    For every $k\in\Z^n$, 
    \[\widehat{\psi_f}(k):=\int_{\T^n}z^{-k}\psi_f(z)d\mathrm{Haar}(z).\]
    The uniform convergence of $\psi_{f_R}$ to $\psi_f$ allows to take the limit $R\to\infty$ under the integral, which gives
    \begin{equation}\label{eq:toruspsicoeff}
        \widehat{\psi_f}(k)=
        \lim_{R\to\infty}b_k f_R(L^tk)=
        b_{k}f(L^tk).
    \end{equation}

    Define the closed set
    \[
        M_0:=\{z\in\T^n:0\in\Lambda_z\}.
    \]
    We use the local point-selection construction
    from the proof of
    Proposition~\ref{prop:lighthousesmooth}.
    Fix $z_0\in M_0$. Common separation and vague
    continuity imply that, for $z$ sufficiently close to $z_0$,
    there is a unique point in the $\delta/4$ ball around the origin, $y(z)\in\Lambda_z\cap B_{\frac{\delta}{4}}$, and no other points in the $\delta/2$ ball,  $\Lambda_z\cap B_{\frac{\delta}{2}}=\{y(z)\}$. Choose $f_1,\ldots,f_d\in\mathcal S(\R^d)$
    whose Fourier transforms are real, such that $\supp{(\widehat{f_j})}\subset B_{\frac{\delta}{2}}$, and $\widehat{f_j}(\xi)=\xi_j$ for all $\xi\in B_{\frac{\delta}{4}}$. Then
    \[
        \psi(z):=(\psi_1(z),\ldots,\psi_d(z)),
        \qquad
        \psi_i(z):=\psi_{f_i}(z)
        =\sum_{x\in\Lambda_z}\widehat f_i(x),
    \]
    satisfies $\psi(z)=y(z)$ near $z_0$.

    Choose $P\in\R^{n\times c}$ whose columns form
    a basis of $\ker(L^t)$ and put
    $z(x)=z_0\exp(2\pi iPx)$ for $x\in\R^c$.
    For $x$ sufficiently small, for all $y\in B_{\frac{\delta}{4}}$,
    \eqref{eq:torusfamilytranslation} gives
    \[
        z(x)\varphi(y)\in M_0
        \quad\Longleftrightarrow\quad
        y\in\Lambda_{z(x)}
        \quad\Longleftrightarrow\quad
        y=\psi(z(x)).
    \]
    Thus, as in Proposition~\ref{prop:lighthousesmooth},
    \[
        x\longmapsto
        z(x)\exp\bigl(2\pi iL\psi(z(x))\bigr)
    \]
    parametrizes $M_0$ locally as a graph transverse
    to $L$. Here continuity of $\psi$ is sufficient:
    no lighthouse assumption or decay of Fourier coefficients is needed. Hence $M_0$ is
    a topological submanifold of dimension $c=n-d$.

    We now identify $M_0$ with 
    $$M:=\overline{\varphi(\Lambda)}.$$ Equation
    \eqref{eq:torusfamilytranslation} gives
    $\varphi(\Lambda)\subset M_0$, and therefore
    $M\subset M_0$. Conversely, let $z\in M_0$
    and choose $\varphi(y_j)\to z$. Vague convergence
    of $\mu_{\Lambda-y_j}$ to $\mu_z$ supplies
    points $x_j\in\Lambda-y_j$ with $x_j\to0$.
    Then $y_j+x_j\in\Lambda$ and
    $\varphi(y_j+x_j)\to z$, proving $M_0\subset M$.
    Consequently, for $M_z=z^{-1}M$,
    \[
        \Lambda_z
        =
        \{x\in\R^d:z\varphi(x)\in M\}
        =
        \Lambda(M_z,L).
    \]
    Taking $z=\mathbf{1}$ proves part~(1).

    For proving parts (3) and (2), let $\nu=\nu_L$ be the $L$-directed
    measure on $M$. We first prove Part (3), namely, showing that $\widehat{\nu}(k)=b_k$.

    To see why, choose the product
    neighborhood from
    Lemma~\ref{lem : product structure}, with the param\-etrization $M\times (-\epsilon_0,\epsilon_0)^d\ni(w,x)\mapsto w\varphi(x) $, so that $d\mathrm{Haar}(w\varphi(x))=d\nu(w)dx$. If $f\in\mathcal S(\R^d)$ with $\widehat{f}$ supported in $(-\epsilon_0,\epsilon_0)^d$, then, $\supp{(\psi_f)}$ is contained in the product neighborhood of $M$, and   
    \[(\psi_f)(w\varphi(x))=\sum_{y\in\Lambda_{w}-x}\hat{f}(y)=\hat{f}(-x),\qquad \forall (w,x)\in M\times (-\epsilon_0,\epsilon_0)^d.\]
Therefore, we calculate $\widehat{\psi_f}$ using this parametrization
    \begin{align*}
        \widehat{\psi_f}(k)
        &=
        \int_{(w,x)\in M\times(-\epsilon_0,\epsilon_0)^d}
        \widehat{f}(-x)
        \bigl(w\varphi(x)\bigr)^{-k}
        \,dx\,d\nu(w)\\
        &= \int_{w\in M}w^{-k}d\nu(w)\int_{x\in (-\epsilon_0,\epsilon_0)^d}
        \widehat{f}(x)e^{2\pi i\langle Lx,k\rangle }
        dx\\
        & = \widehat\nu(k)f(L^tk).
    \end{align*}
   Comparing with \eqref{eq:toruspsicoeff}, and
    for each $k$ choosing such an $f$ with
    $f(L^tk)\ne0$, gives the needed result
    $\widehat\nu(k)=b_{k}$.

Finally, $\mu_\Lambda$ is Bohr almost periodic by
Corollary~\ref{cor:crystallineBAP}. Its Bohr--Fourier
coefficients coincide with the Fourier coefficients of its
distributional Fourier transform
\cite{meyer2018global}*{Lemma 5.11}. Hence, for every
$k\in\Z^n$,
\[
    \lim_{R\to\infty}
    \frac{1}{\mathrm{vol}(B_R)}
    \sum_{x\in\Lambda\cap B_R}
    e^{-2\pi i\langle x,L^tk\rangle}
    =
    \widehat{\mu_\Lambda}(\{L^tk\})
    =
    b_k
    =
    \widehat\nu(k).
\]
Thus the averaged-pushforward measures converge to
$\nu$ against every character. The case $k=0$ gives
convergence of their total masses, so positivity
and the Stone--Weierstrass theorem imply weak
convergence to $\nu$, proving part~(2).
\end{proof}

\bibliographystyle{plain}
\bibliography{biblio}
	
\end{document}

%% file: operators.tex
\usepackage{mathrsfs}

\newcommand{\Hom}{\operatorname{Hom}} 

\newcommand{\rank}{\operatorname{rank}} 

\newcommand{\Gr}{\operatorname{Gr}}

\newcommand{\T}{\operatorname{T}}

\newcommand{\tr}{\operatorname{tr}}



%% file: fonts.tex
\newcommand{\cA}{{\mathcal A}}
\newcommand{\cB}{{\mathcal B}}

\newcommand{\cS}{{\mathcal S}}

\newcommand{\fp}{{\mathfrak p}}

\newcommand{\C}{{\mathbb C}}
\newcommand{\R}{{\mathbb R}}

\newcommand{\Q}{{\mathbb Q}}
\newcommand{\N}{{\mathbb N}}
\newcommand{\Z}{{\mathbb Z}}

